\documentclass[11pt,a4paper]{article}

\usepackage[british]{babel}
\usepackage[T1]{fontenc}
\usepackage[utf8]{inputenc}
\usepackage{microtype}            % better kerning & protrusion

\usepackage{newpxtext}            % Palatino-like text
\usepackage{newpxmath}            % and matching maths
\usepackage[a4paper,margin=28mm]{geometry}
\usepackage[skip=6pt plus 2pt, indent=0pt]{parskip}

\usepackage{amsmath,mathtools}
\numberwithin{equation}{section}  % optional; remove if you prefer global numbering

\usepackage{graphicx}
\usepackage{wrapfig}
\usepackage{float}
\graphicspath{{figs/}}            % put your images in ./figs/
\usepackage{booktabs}
\usepackage[labelfont=bf,textfont=it]{caption}
\usepackage[dvipsnames]{xcolor} 
\usepackage[colorlinks=true, allcolors=blue!50!black]{hyperref}
\usepackage[nameinlink,capitalise,noabbrev]{cleveref} % \Cref, \cref for “Figure 1”, etc.

\usepackage[numbers,sort&compress]{natbib}
\title{Crocheting Mathematics}
\author{Hanne Kekkonen}
\date{\small Preprint of a chapter published in \textit{Shapes in Action}, Aalto ARTS Books, 2023.}

\begin{document}

\maketitle

\vspace{0.5cm}

\begin{abstract}
\noindent Crochet provides a superior method for the production of two-dimensional surfaces from one-dimensional material. Compared to any of the other known processes to generate constant flat, spherical or hyperbolic shapes, it is the most flexible and precise way to build a dynamical system with very simple local rules and with very high precision of the intrinsic curvature.
\end{abstract}

%%%%%%%%%%%%%%%%%%%%%%%%%%%%%%%%%%%%%%%%%%%
\section{Short history of Mathematical Crochet}

Crocheting is an inherently mathematical process. You need to keep track of the number of stitches added to a disc to make a tablecloth or a hat. The first should have zero curvature, whereas the second should have the correct positive curvature. Complex lace crocheting instructions can look like something you would need a university degree to decipher, and to understand a written pattern you must learn a new language. Crocheting also allows the creation of many complex mathematical objects that are very difficult to make in any other way, as long as you understand the underlying properties of those shapes.

Perhaps the most famous crocheted mathematical model was created by Daina Taimina in 1997, when she realised that one could crochet a hyperbolic plane instead of making a rather fragile and laborious paper-strip model. The formula for the model is based on hyperbolic geometry, more specifically on the circumference of circles on hyperbolic surfaces. The paper describing how to crochet a hyperbolic surface was published in The Mathematical Intelligencer in 2001 \cite{Henderson2001}, and, as the authors Daina Taimina and David Henderson predicted, it inspired many other crocheting mathematicians to create new mathematical models. In 2004, a second crocheted model appeared on the cover of The Mathematical Intelligencer. This time, Hinke Osinga and Bernd Krauskopf described how to crochet an approximation of the Lorenz manifold \cite{Osinga2004}.

The idea of using crochet or knitting to produce mathematical models is, however, much older. The Scottish chemist Alexander Crum Brown knitted several interlinked surfaces to visualise, and perhaps further, his own understanding of the ideas presented in his paper \textit{On a Case of Interlacing Surfaces} in the late nineteenth century \cite{Brown1886, Dunning2015}. The British mathematician Miles Reid wrote a paper on knitting mathematical surfaces in the 1970s \cite{Reid1971}, which has since inspired several new patterns, including a Möbius scarf and different ways to knit the Klein bottle. There are also many more recent crocheted and knitted mathematical models, including the real projective plane by Claire Irving \cite{Irving2004}, several topological sculptures by Shiying Dong \cite{Dong2023}, and a bitruncated pentachoron, a projection of a four-dimensional polytope into our three-dimensional space, by Kirsi Peltonen \cite{Luotoniemi2019}. See also \cite{Taimina2018} for further examples.

In this chapter, we use crochet to create different types of mathematical models. Even though crochet and knitting may seem rather similar, there are several differences between the two. The main one is that in crochet you work with only one `live' stitch; the previous stitch is always completed before moving to the next one. In knitting, all the stitches in a round stay on the needles, and you must bind them off separately if they are not needed on the next round. Also, all knitting stitches are the same height, whereas there are several different crochet stitches, which allows part of a round or row to have a different height. These differences make crocheting much more agile when creating complicated shapes. The fabric produced by knitting or crocheting is also quite different: knitting creates a supple, stretchy fabric that is well suited for clothing, while crocheted fabric has little stretch and can be made quite rigid using a small hook. Using a special crochet stitch discussed in the next section, the surface can be made smooth and very stiff, allowing even larger structures to keep their shape without supporting wire frames.

%%%%%%%%%%%%%%%%%%%%%%%%%%%%%%%%%%%%%%%

\section{Surfaces with only one side}

%\begin{wrapfigure}{r}{5cm}
%\includegraphics[width=5cm]{figures/ClocksColumn}
%\caption{You cannot define clockwise consistently on the M\"obius band because it is non-orientable.}
%\label{Fig:Clocks}
%\end{wrapfigure} 

Take a strip of paper, give it a half twist, and glue the ends together. You have just created a shape called the Möbius band, a very simple-looking object with rather perplexing properties. The Möbius band is an example of a one-sided surface, and it was first described in 1858 by the German mathematician and astronomer August Möbius. The shape was actually discovered independently by another German mathematician, Johann Listing, a bit earlier, but he did not publish his work until 1861. The history of the shape goes back much further than the 1850s, and it can already be found in old Roman mosaics. It is often used to symbolise infinity and endless cycles, such as the passing of the seasons or, slightly less poetically, recycling. The Möbius band has also been considered in mechanics and engineering. The first known example appears in the \textit{Book of Knowledge of Ingenious Mechanical Devices} by al-Jazari in 1206 \cite{Ibn1206}, in which a chain pump with the chain arranged as a Möbius band is shown in several illustrations. Drive and conveyor belts with Möbius-band topology were also described in the late 19\textsuperscript{th} and early 20\textsuperscript{th} centuries see, e.g., \cite{Inman1903}. The advantage of such belts is that wear is distributed over `both sides'. See, e.g., \cite{Cartwright2016} for more on the historical use of the Möbius band.

\begin{figure}[H]
\centering
\begin{minipage}[b]{0.95\textwidth} 
	\includegraphics[width=\textwidth]{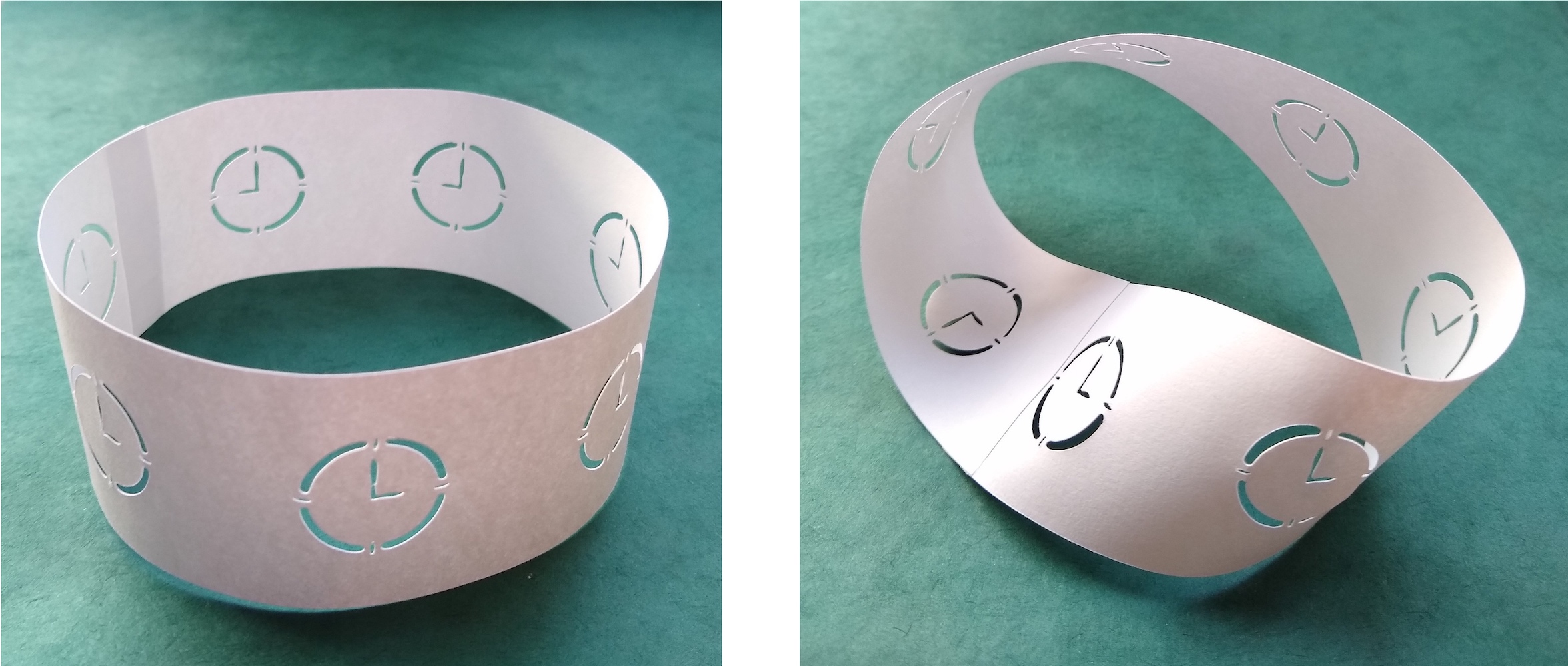}
\end{minipage}
\caption{You cannot define clockwise consistently on the Möbius band (right) because it is non-orientable.}
\label{Fig:MobiusClocksRow}
\end{figure}

To see that the loop you just created really has only one side, you can use a pencil to draw a line along the centre of the strip. You will see that the line runs through the whole loop; there is no separate inside and outside. This one-sidedness in our Euclidean 3-space makes the Möbius band a non-orientable surface. An orientable surface allows a consistent orientation to be assigned over the entire surface. To better understand the concept, we can consider the difference between a cylinder and the Möbius band.

We begin by placing a clock face, fixed at 3 o’clock, on the surface of a cylinder, as shown on the left in Figure \ref{Fig:MobiusClocksRow}. If you slide the clock along any closed curve on the surface, keeping the hour hand pointing along the curve, then the clock face continues to show 3 o’clock. If you move the minute hand to five past, it moves clockwise on all the clocks, meaning that we can define ‘clockwise’ and ‘anticlockwise’ consistently.

If you place the same clock face on a Möbius band and slide it around the central curve, you will notice that when the clock first returns to its starting point, it has turned upside down. If you moved the minute hand clockwise in the first position, it would move anticlockwise in the final clock. This means that we cannot consistently define ‘clockwise’ and ‘anticlockwise’, and so the Möbius band is a non-orientable surface.

\begin{wrapfigure}{l}{5.4cm}
\includegraphics[width=5.4cm]{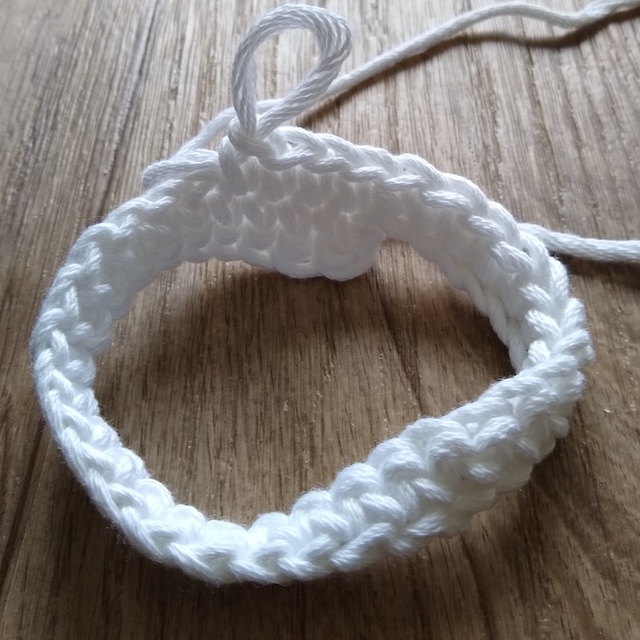}
\caption{You can crochet a Möbius band by crocheting a row, giving it a half twist, and connecting the ends into a loop. You can then finish the first round by crocheting on the ‘wrong’ side through the original chain.}
\label{Fig:CreateMobius}
\end{wrapfigure} 

You can create a more durable Möbius band by crocheting. One way to do this is to follow the same construction as with the paper model: crochet a wide strip, give it a half twist, and then stitch the ends together. Another way is to exploit the fact that the Möbius band has only one boundary: start by crocheting one row, give it a half twist, and then connect the ends to create a loop as shown in Figure \ref{Fig:CreateMobius}. You can then simply crochet along the single edge until the desired width is achieved. This second method also allows you to crochet other versions of the Möbius band, as we will soon see.

Most materials have two distinctively different sides. You are most probably wearing a knitted shirt. If you look very carefully at the right side of your shirt, you will see small V shapes. If you then look at the wrong side, you will notice that instead of little Vs there are small bumps. This is also true for crocheted surfaces.\footnote{There are no crocheting machines, so if you ever receive a crocheted gift you can be sure it is handmade.}

\begin{figure}[t]
\centering
\begin{minipage}[b]{0.93\textwidth} 
\includegraphics[width=\textwidth]{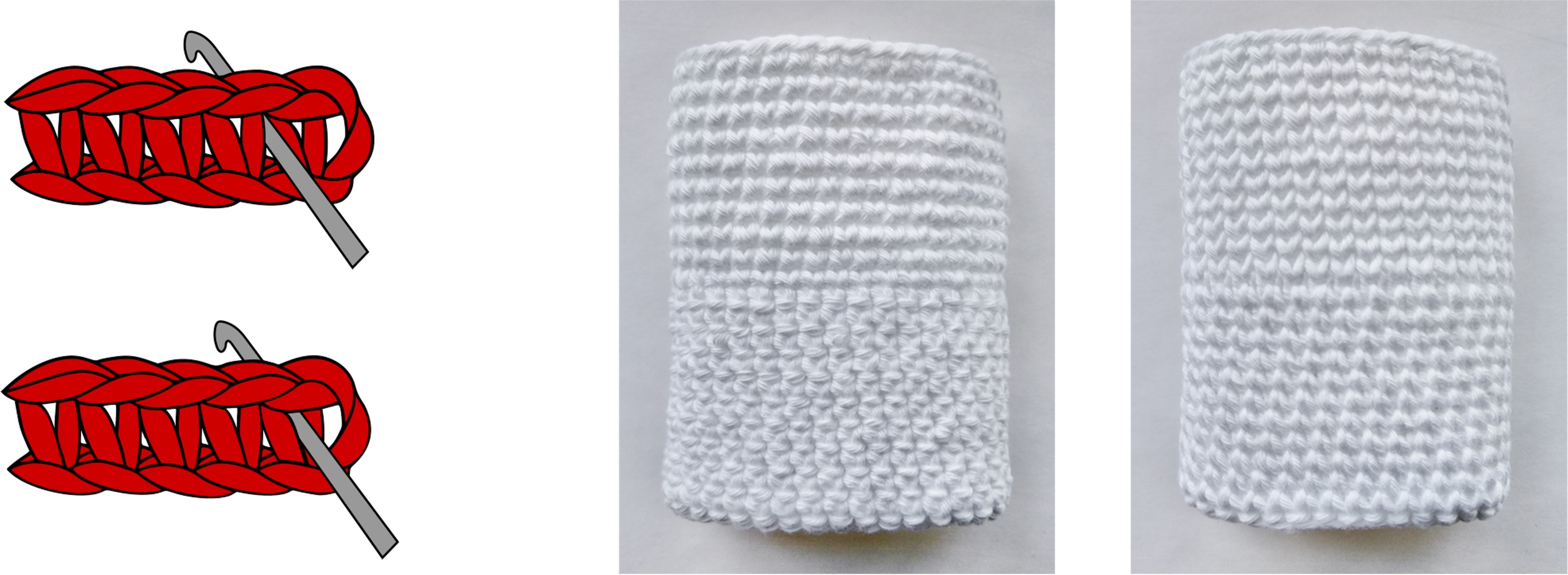}
\end{minipage}
\caption{In the traditional single crochet stitch, the yarn is pulled under the loop connecting two V-shaped stitches (top left). This creates a knobbly surface with a clearly different front (top of the left tube) and back (bottom of the left tube). In the split stitch, the yarn is pulled through a V-shaped stitch (bottom left). The resulting surface is smoother and more rigid, with less distinct sides (right).}\label{Fig:HyperbolicFootball}
\label{Fig:Stitches}
\end{figure}

The traditional crochet stitch is created by pulling the yarn through the loop connecting two stitches on the previous row, illustrated in the top left in Figure \ref{Fig:Stitches}. This creates a rather bumpy surface, with the right side showing small unstacked V shapes, whereas the wrong side features additional short horizontal lines. Models created with the traditional crochet stitch tend to be quite floppy and need supporting wire to keep their shape. To create more rigid models with similar sides, you can use the split stitch (also known as the waistcoat stitch). This stitch is created by pulling the yarn through the centre of the stitch on the previous row. The resulting surface looks like knitted fabric on the right side, and the wrong side does not have horizontal lines, making the two sides quite similar. The surface is also very stiff, so your models will not collapse under their own weight.

 \begin{wrapfigure}{r}{5cm}
\includegraphics[width=5cm]{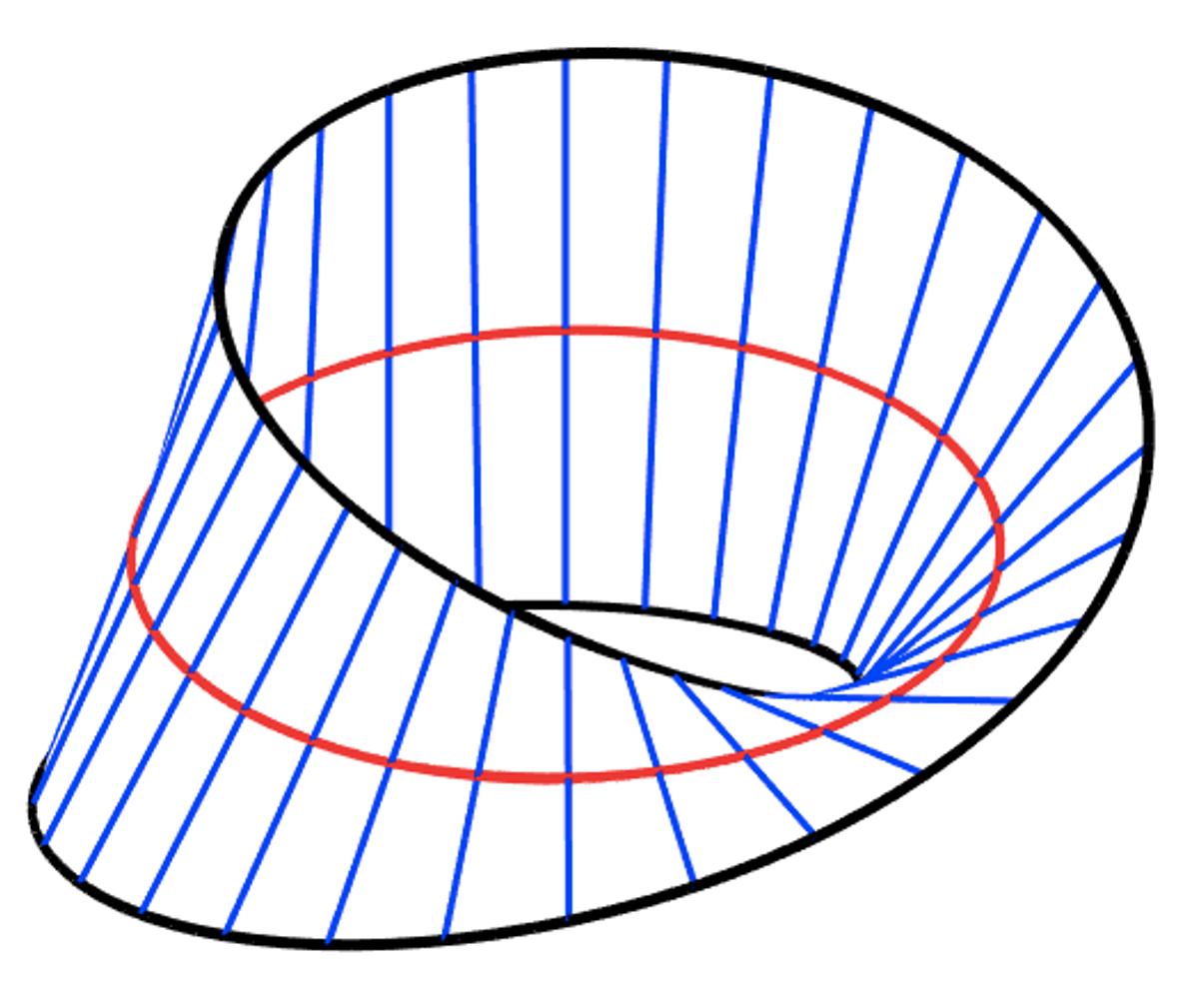}
\caption{The parametrisable M\"obius band.}
\vspace{-0.9cm}
\label{Fig:Mobius2}
\end{wrapfigure} 

Interestingly, even though it is very easy to make a Möbius band from a strip of paper, there is no simple parametrisation for that model. However, there is another embedding of the Möbius band that can be described with a simple set of equations:
\begin{align*}
x(r, \theta) & =\left(1+r \cos(\theta)\right) \cos (2\theta) \\
y(r, \theta) & =\left(1+r \cos(\theta)\right) \sin (2\theta) \\
z(r, \theta) & =r \sin (\theta)
\end{align*}
where $\theta \in[0,\pi]$ and $r \in[-1, 1]$, as shown in Figure \ref{Fig:Mobius2}. This model looks rather similar to the one we considered earlier, but you cannot make it from a flat strip of paper. If you cut this Möbius band open, it would not lie flat like the paper model. Also, its centre forms a perfect circle, and you can draw straight lines through this centre circle, making the embedding a \textit{ruled surface}.

\begin{figure}[t]
\centering
\begin{minipage}[b]{0.95\textwidth} 
	\includegraphics[width=\textwidth]{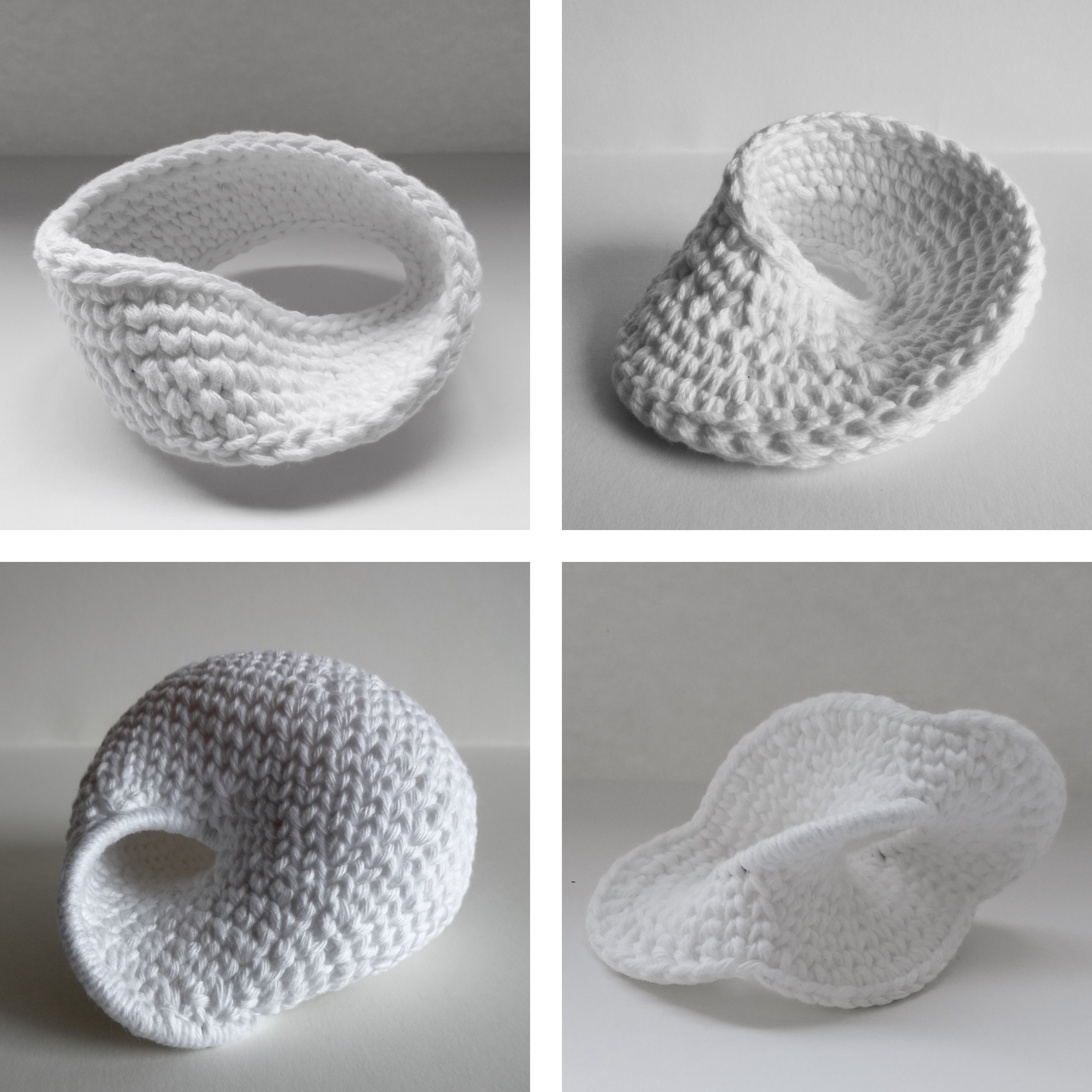}
\end{minipage}
\caption{Crocheted `paper strip' Möbius band (top left), the parametrisable M\"obius band (top right), Möbius snail (bottom left), and a Möbius band with a circular boundary at the centre and a planar-like end stretching away from it. The outer edge of the planar Möbius band should extend to infinity, so the only actual boundary of the shape is the central circle.}\label{Fig:CrochetedMobiusBands}
\end{figure}

Crocheting the parametrised Möbius band requires a bit more deliberation. The first method described earlier is out of the question, but the second can be modified to approximate this version of the Möbius band. The difference from the paper-strip model is that you now have to add and subtract stitches to create the correct shape. Because the surface has varying curvature, we need to calculate the exact positions for the increases and decreases. We discuss calculating exact crochet shapes in more detail in Section \ref{Sec:CalculatingExactCrochetPatterns}.
In the top row of Figure \ref{Fig:CrochetedMobiusBands}, you can see a crocheted version of both embeddings of the Möbius band.

%\begin{figure}[H]
%\centering
%\begin{minipage}[b]{0.95\textwidth} 
%	\includegraphics[width=\textwidth]{figures/MobiusBands}
%\end{minipage}
%\caption{Crocheted `paper strip' Möbius band (left) and the parametrisable M\"obius band (right).}\label{Fig:HyperbolicFootball}
%\label{Fig:MobiusBands}
%\end{figure}
% 
%\begin{figure}[H]
%\centering
%\begin{minipage}[b]{0.95\textwidth} 
%	\includegraphics[width=\textwidth]{figures/MobiusRound}
%\end{minipage}
%\caption{Crocheted Möbius snail (left) and a Möbius band with a circular boundary at the centre and a planar-like end stretching away from it. The outer edge of the planar Möbius band should extend to infinity, so the only actual boundary of the shape is the central circle.}
%\label{Fig:MobiusRound}
%\end{figure}

As mentioned above, the Möbius band has only one boundary component. The boundary is topologically equivalent to a circle, which you can observe by cutting off a very thin slice along the boundary of your paper Möbius band. This means that it is possible to embed the Möbius band in three dimensions so that its boundary forms a perfect circle. One example of such a shape is the Möbius snail, or Sudanese Möbius band, named after the topologists Sue Goodman and Daniel Asimov, who discovered the shape in the 1970s, as shown in the bottom left in Figure \ref{Fig:CrochetedMobiusBands}. This ‘shell’ has two openings that smoothly connect the ‘inside’ and ‘outside’. If you carefully follow the surface, you will notice that, like any Möbius band, it indeed has only one side. You do not need to go over an edge to move from ‘inside’ to ‘outside’, as you would with a traditional shell.

Imagine making a tiny puncture at the bottom of the Möbius snail and stretching the hole bigger and bigger until it extends to infinity. This would create a one-sided shape, a ‘Möbius plane’, with a circular boundary at the centre and the surface stretching to infinity like a plane, as shown in the bottom right in Figure \ref{Fig:CrochetedMobiusBands}. Since we made a hole in the Möbius snail, this new shape is not topologically equivalent to the Möbius bands we have considered so far. It is closely related to the Möbius snail, however, as both are stereographic projections of the same Möbius band on the three-dimensional sphere in four-dimensional space into our three-dimensional space. For the Möbius snail, the projection point is chosen so that it does not lie on the embedded Möbius strip, whereas for the ‘Möbius plane’ the projection point is chosen to produce the most symmetrical image, with the point lying on the strip mapped to infinity.

\begin{wrapfigure}{r}{6cm}
\includegraphics[width=6cm]{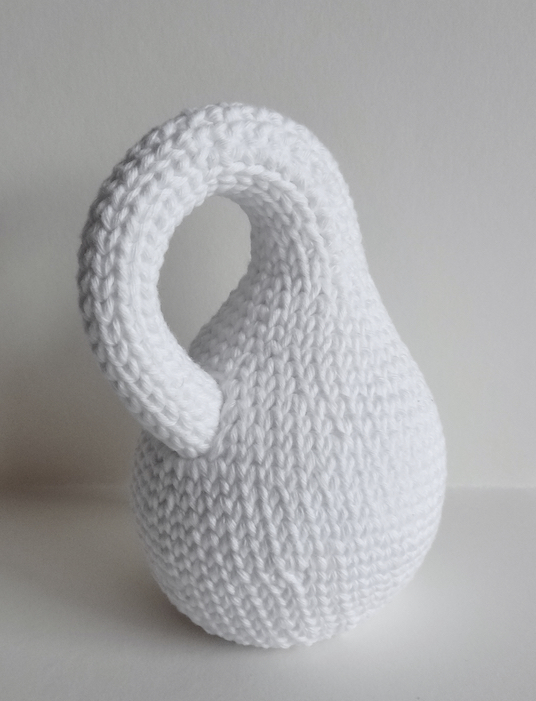}
\caption{A Crocheted Klein bottle.}
\label{Fig:KleinBottle}
\end{wrapfigure} 

If you take two Möbius strips, one with a left-handed half-twist and the other with a right-handed half-twist, and glue them together, you obtain a one-sided surface without a boundary. However, to perform this gluing properly you must move to four-dimensional space, since in our three-dimensional space the procedure produces a shape with a self-intersection. 
The surface, both in four dimensions and the self-intersecting three-dimensional immersion, is known as the Klein bottle, shown in Figure \ref{Fig:KleinBottle}. Even though the Klein bottle does look somewhat like a peculiar bottle, the name is often said to be the result of a translation error, mixing the German word Fläche (surface) with Flasche (bottle).

The concept of the Klein bottle was first detailed by the German mathematician Felix Klein in 1882, who visualised it by inverting a piece of rubber tubing and letting it pass through itself so that the outside and inside meet.
Any small patch of the Klein bottle is two-dimensional, making it an example of a two-dimensional manifold that can only exist in four dimensions. Like the Möbius strip, the Klein bottle is also non-orientable.

%%%%%%%%%%%%%%%%%%%%%%%%%%%%%%%%%%%%%%%

\section{Knots and Seifert Surfaces}

The wheel or the use of fire are often considered to be the greatest inventions of humankind, but we would not have got far without knots. Knots allowed sharp stone tools to be tied to handles making them far more efficient, they enabled the creation of fishing nets and snares that provided larger catches, the building of rafts that allowed early seafaring, and the construction of efficient shelters that enabled humans to live almost anywhere on the globe. The list goes on, and almost any early invention included knots in some way. Since knots were tied into pliable organic materials, which are highly perishable, we have only indirect archaeological evidence of their early use. For example, perforated beads and pendants, the oldest dating to $300\ 000$ BCE, were likely strung on some form of tied cord. The earliest direct evidence of the use of strings and knots comes from several carved human figurines from $22\ 000$–$28\ 000$ years ago wearing elaborate belts and straps, and the oldest fragments of rope are from about $15\ 000$–$17\ 000$ BCE.

Besides their primary function of tying objects together, knots have also been used to record information and for decoration. The Incas used quipu, string- and knot-based devices, for everything from keeping records and calendrical information to monitoring tax obligations and military organisation. The strings stored information encoded as knots, and the largest quipu could have thousands of cords. Similar devices were also used in ancient China and Japan.

Knots have been used as decoration throughout history. For example, many Babylonian and Assyrian carvings depict clothing with elaborate knot adornments. Knotted fringes were also common among Arab weavers and could be found on towels, shawls, and veils. This type of decorative knot is nowadays often referred to as macramé. For more on the history and use of knots, see, e.g., \cite{Turner1996}.

%\begin{wrapfigure}{r}{7cm}
%\includegraphics[width=7cm]{figures/KnotTable}
%\caption{All prime knots with seven or less crossings.}
%\label{Fig:KnotTable}
%\end{wrapfigure} 

\begin{figure}[H]
\centering
\begin{minipage}[b]{0.95\textwidth} 
	\includegraphics[width=\textwidth]{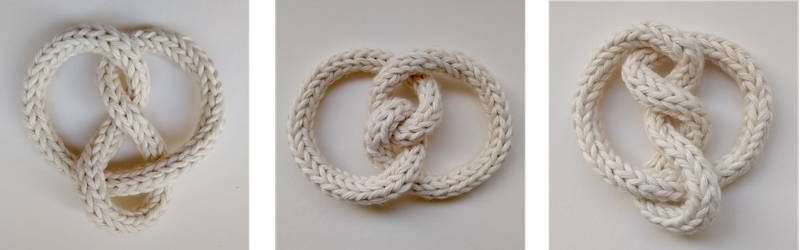}
\end{minipage}
\caption{Three different projections of the figure eight or $4_1$ knot. Notice how closely the rightmost version looks like the $5_2$ knot in Figure 8. The tube is made using a spool knitting tool.}\label{Fig:FigureEightKnots}
\end{figure}

Even though humans have used knots since the dawn of time, the mathematical study of knots began only in the early 19\textsuperscript{th} century with the work of the German mathematician Carl Friedrich Gauss. The real push for modern knot theory came from the suggestion by the Scottish mathematician–physicist Lord Kelvin that atoms might be knots in the aether. This theory, along with the aether, was later discarded, but it led another Scottish mathematician–physicist, Peter Guthrie Tait, to create the first systematic classification table of knots. For the first hundred years, knot theory was considered to be of mainly mathematical interest, but in the late 20\textsuperscript{th} century it became connected to quantum field theory and hyperbolic geometry \cite{Thurston2022}, and found applications in biology, chemistry, and mathematical physics \cite{Buck2009, Horner2016}. Knot theory has even been used to study the mathematical properties of knitted materials \cite{Singal2023}.

\begin{figure}[H]
\centering
\begin{minipage}[b]{0.9\textwidth} 
	\includegraphics[width=\textwidth]{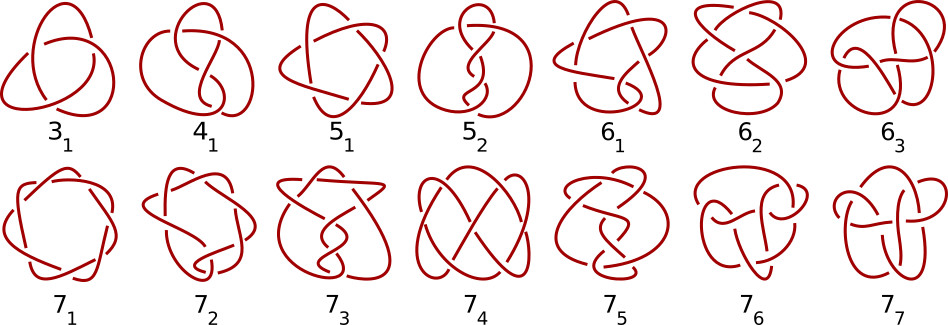}
\end{minipage}
\caption{All prime knots with seven or less crossings.}
\label{Fig:KnotTable}
\end{figure}

But what exactly do mathematicians mean by a knot? The main difference between a mathematical knot and a standard knot is that mathematical knots are closed. The simplest knot is just a circle, called the unknot or trivial knot. One of the main questions in knot theory is whether a closed curve is truly knotted or can be untangled, in other words, whether you can deform the curve in space into the unknot without breaking the loop. More generally, we can ask whether two given curves represent different knots or are actually the same knot, in the sense that one can be homeomorphically deformed into the other. In Figure \ref{Fig:FigureEightKnots} you can see three different ways to arrange the $4_1$ knot, often called the figure-eight knot. Another question in knot theory considers how many different knots exist under given constraints. Knots can be classified by the number of times the string crosses over itself. The knots in Figure \ref{Fig:KnotTable} are called prime knots because they cannot be formed by joining smaller knots together.

You can easily make mathematical knots by using a technique called spool knitting (also known as French knitting) or by crocheting a so-called i-cord. Once you have made a long enough tube, you can form a knot using the knot table in Figure \ref{Fig:KnotTable} as a guide. To finish the knot, tie or stitch the ends of the tube together. But how long is a `long enough' tube? Unlike mathematical knots, knots created from any physical string have thickness, so there is a lower bound for the required length.

A classical knot theory question asks: ``Can I tie a knot in a foot-long rope that is one inch thick?" This very innocent-sounding question is mathematically rather difficult to solve. To answer it, we need to consider the so-called ropelength, which is given by the length of the rope divided by its radius. The ropelength of a foot-long ($30.48$ cm) and one-inch ($2.54$ cm) thick rope is $30.48 / 1.27 = 24$. Computer simulations have shown that the minimum ropelength required to tie the trefoil knot is very close to $32.74$ \cite{Ashton2011, Rawdon2003, Stasiak1998}. The sharpest theoretical lower bound for the ropelength of any non-trivial knot is $31.32$ \cite{Denne2006}, which is very close to the bound given by simulations. So the answer to the question is no: a foot-long, one-inch thick rope is far too short to tie a knot.

The thickness of the tube created using a four-pin knitting spool and aran-thickness cotton yarn is about $0.8$ cm, so the minimum length of your tube should be $32.8 \times 0.4\ \text{cm} \approx 13$ cm. This would create a very tight knot, so I recommend doubling the minimum length for the trefoil knot and adding about $10$ cm for every extra crossing. If you would like to make a set of tight knots, you can find a list of computer-simulated minimum ropelength bounds for knots up to $10$ crossings in \cite{Ashton11b}.

%\begin{wrapfigure}{r}{4cm}
%\includegraphics[width=4cm]{figures/CreatingKnot}
%\caption{You can crochet a  Seifert surface of the trefoil knot by making two small discs and connecting them by three `legs' that all have half twist to the same direction (top). You can make the shape smoother by using double and treble stitches close to the joints on the next round (bottom).}
%\label{Fig:CreatingKnot}
%\end{wrapfigure}  

\begin{figure}[H]
\centering
\begin{minipage}[b]{0.97\textwidth} 
	\includegraphics[width=\textwidth]{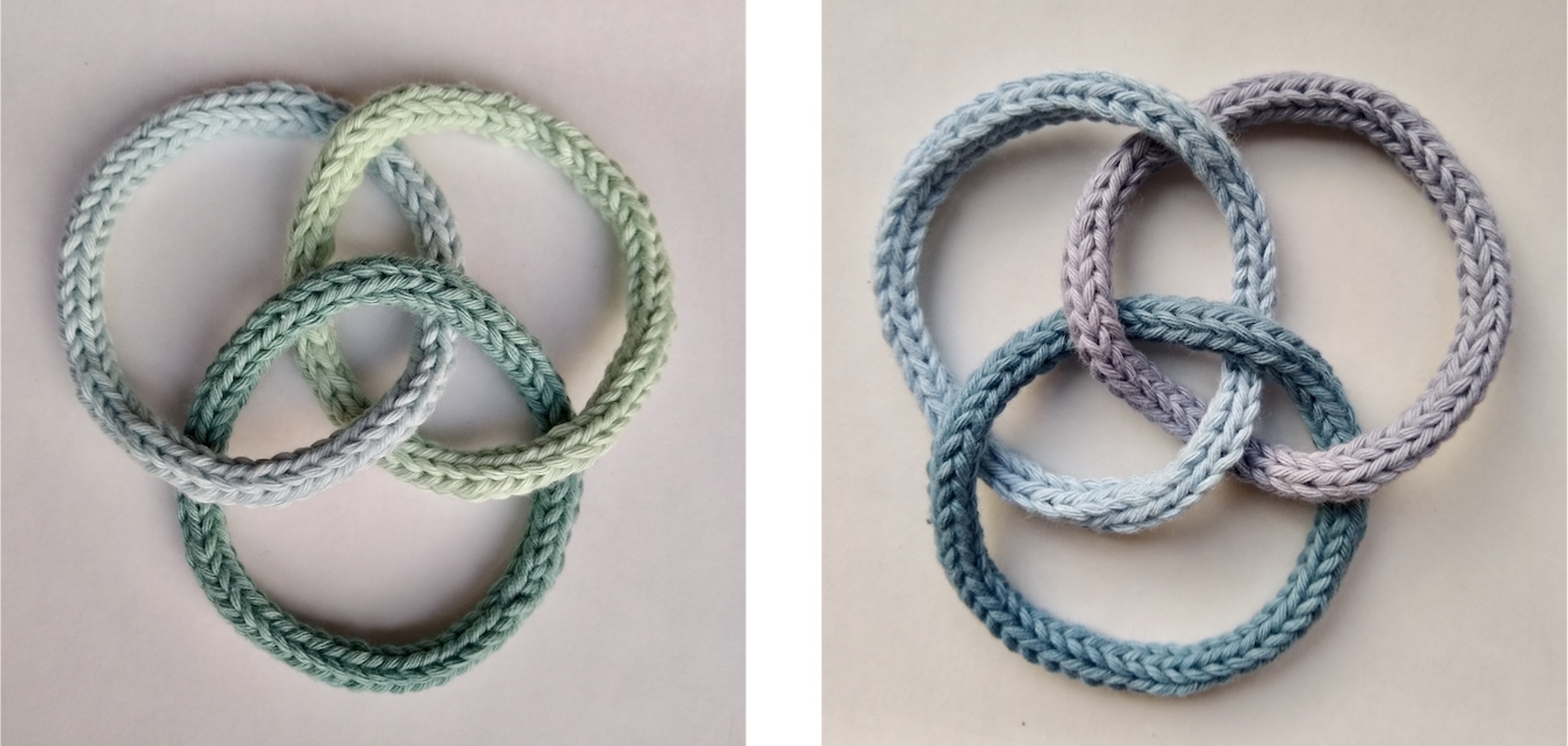}
\end{minipage}
\caption{The $(3,3)$ torus link and the Borromean rings.}\label{Fig:Links}
\end{figure}

Another closely related mathematical concept is that of a link, which is a collection of knots that do not intersect but may be linked together. The simplest non-trivial link is the Hopf link, which consists of two circles linked together once. The $(3,3)$ torus link and the Borromean rings both consist of three interlocked rings, as shown in Figure \ref{Fig:Links}, but while the three rings in the $(3,3)$ torus link are all pairwise linked, no two of the Borromean rings are linked with each other. This means that if you remove one of the rings, the other two will fall apart. Because of this, the Borromean rings have been used around the world to symbolise strength in unity.

\begin{figure}[H]
\centering
\begin{minipage}[b]{0.95\textwidth} 
	\includegraphics[width=\textwidth]{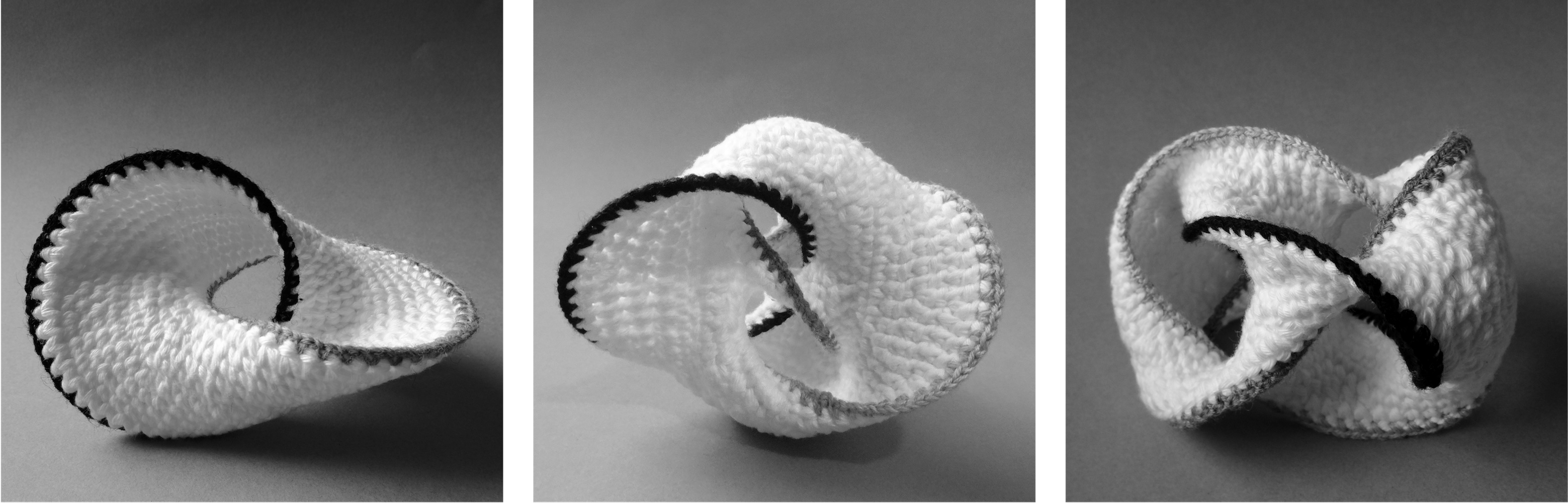}
\end{minipage}
\caption{Seifert surfaces of the $(2,2)$ torus link or Hopf link (left), $(3,3)$ torus link (middle), and the Borromean rings (right).}
\label{Fig:SeifertLinks}
\end{figure}

The Borromean rings are actually the simplest example of a Brunnian link, in which the link becomes a set of unlinked loops if any one component is removed. Another difference between the two is that you can create the $(3,3)$ torus link from perfectly round circles, whereas the components of the Borromean rings can never be perfectly round. Even though the two models are rather different, it can be difficult to tell them apart when the links lie flat. This is also true more generally: differentiating between two knots can be quite difficult just by looking at a flat representation of them. That is why we next consider Seifert surfaces, which can be used to study the properties of the associated knots or links.

A Seifert surface is an orientable surface whose boundary is a given knot or link. This means that the Möbius band is not a Seifert surface, even though its boundary is equivalent to the unknot. Seifert surfaces can be used to study the properties of knots and links, but they are also interesting in their own right. It can be difficult to tell the $(3,3)$ torus link and the Borromean rings apart in Figure \ref{Fig:Links}, but by looking at their Seifert surfaces in Figure \ref{Fig:SeifertLinks} you can see a clear difference. Both surfaces have three different levels joined together, but the Borromean rings have six twisted legs joining the middle part to the top and bottom, whereas the $(3,3)$ torus link has a more complex-looking structure.

\vspace{3mm}
\begin{figure}[H]
\centering
\begin{minipage}[b]{0.97\textwidth} 
	\includegraphics[width=\textwidth]{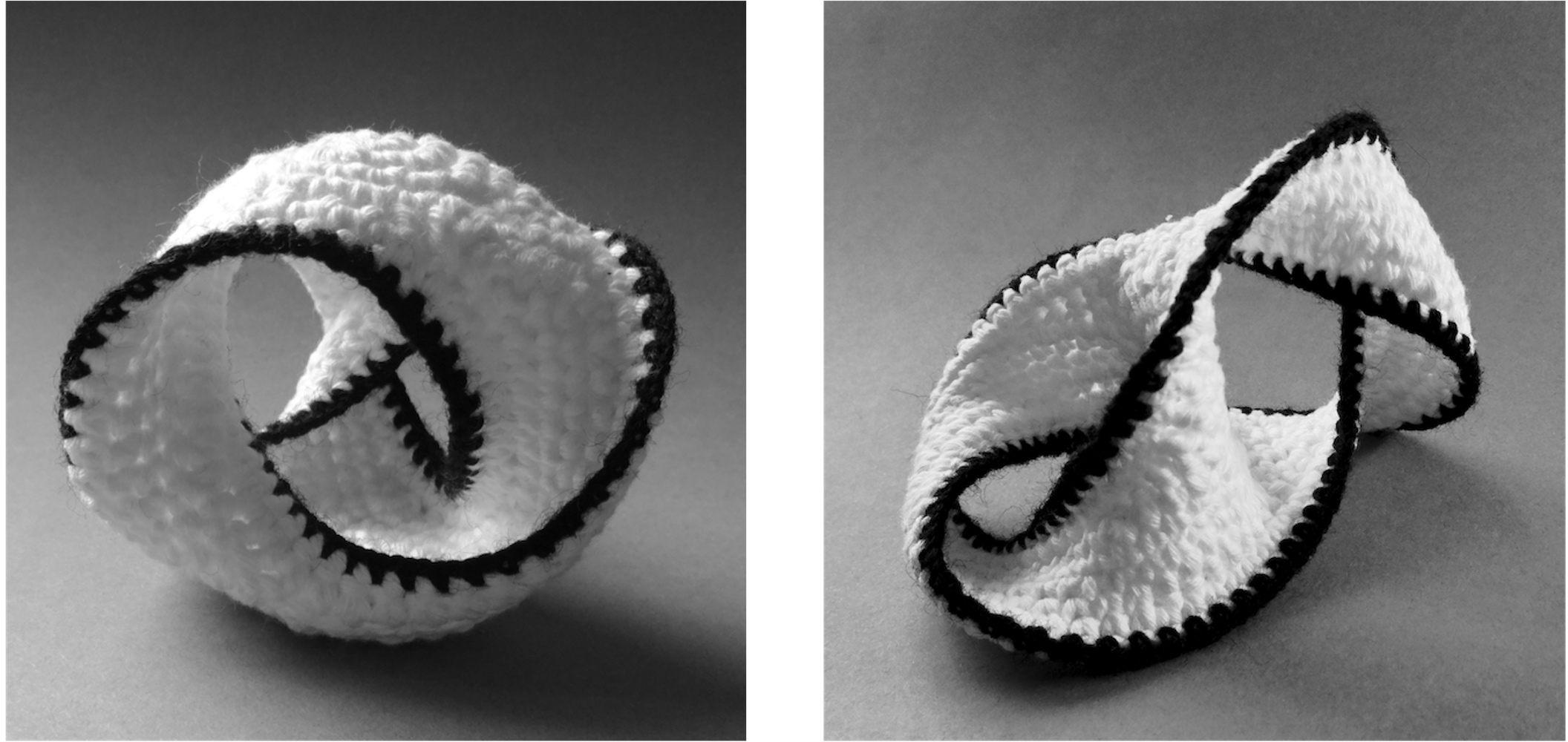}
\end{minipage}
\caption{Seifert surfaces of the trefoil knot (left) and the figure eight knot (right).}
\label{Fig:SeifertKnot}
\end{figure}

Every knot has an associated Seifert surface, and the so-called Seifert algorithm can be used to produce a Seifert surface for a given knot \cite{Seifert1935}. Seifert surfaces are not unique: a single knot or link can have many inequivalent Seifert surfaces. For a more detailed introduction to knot theory, see, e.g., \cite{Adams2004}, which is freely available online.

\begin{figure}[H]
\centering
\begin{minipage}[b]{0.95\textwidth} 
	\includegraphics[width=\textwidth]{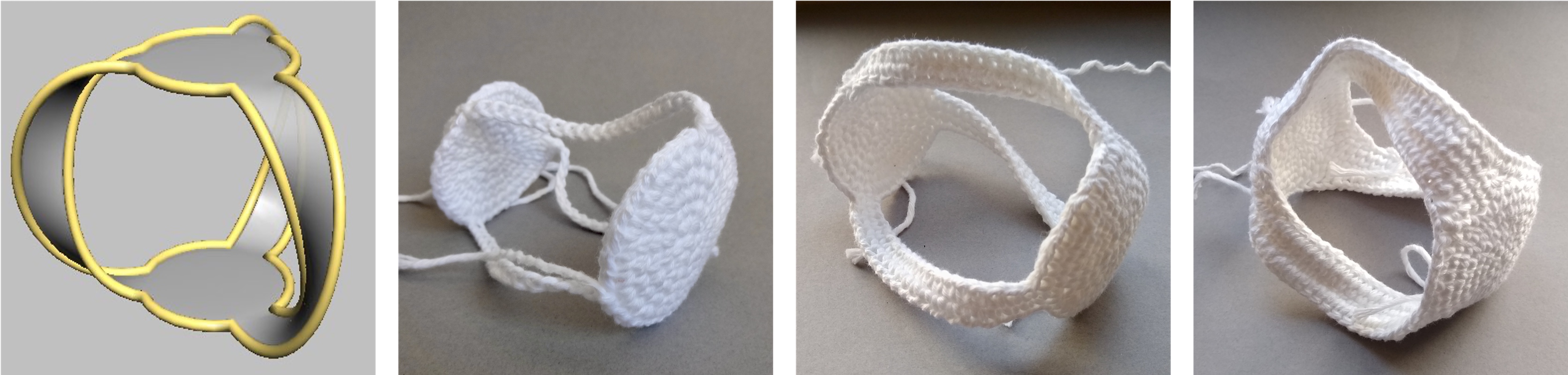}
\end{minipage}
\caption{The SeifertView on the 'crochet instruction' mode (left).
You can crochet a  Seifert surface of the trefoil knot by making two small discs and connecting them by three `legs' that all have half twist to the same direction. The shape can be made smoother by using double and treble stitches close to the joints on the second round (right).}
\label{Fig:CreatingKnot}
\end{figure}

You can crochet your own Seifert surfaces by connecting small discs with twisted chains and crocheting along the boundary until you achieve the desired size, as shown in Figure \ref{Fig:CreatingKnot}. I recommend crocheting a thin plastic cable (e.g., transparent trimmer line) into the model during the last round to give it a cleaner shape, and adding a darker-coloured rim to highlight the knot. Figure \ref{Fig:SeifertKnot} shows the Seifert surfaces of the trefoil knot and the figure-eight knot. If you would like to explore more surfaces associated with knots and links, you can use a program called SeifertView (for Windows) to create your own Seifert surfaces \cite{Wijk2016}.

%%%%%%%%%%%%%%%%%%%%%%%%%%%%%%%%%%%%%%%

\section{Calculating exact crochet patterns}\label{Sec:CalculatingExactCrochetPatterns}

\subsection{Surfaces with constant curvature}

So far, we have mostly discussed models that illustrate a mathematical idea and can be crocheted ‘freestyle’ without a strict pattern. A Klein bottle can be short and wide or tall and narrow, and the Seifert surface of a knot can have flat top and bottom sections or more domed ones. But there are also many shapes that require calculating exact crochet instructions to achieve the intended form. As a simple example, consider a flat disc. To achieve zero Gaussian curvature, we can use the fact that the circumference of a circle on a Euclidean plane is $2\pi R$, where $R$ is the radius of the circle. If the height of your stitch is $H$, the circumference of a crocheted disc on round $\ell$ should be $2\pi (\ell \cdot H)$. To get the number of required stitches, divide this by the width of your stitch $W$; the number of stitches is therefore $2\pi (\ell \cdot H)/W$, rounded to the nearest integer. If the width and height of your stitch are about the same, we can use the approximation $\pi \approx 3$ to deduce the well-known instruction for crocheting a flat disc using the single crochet stitch: start by crocheting $6$ stitches into a ring, and add $6$ stitches on every round. Since a disc looks the same in every direction, you should space the increases evenly.
Notice that the amount of added stitches depends linearly on the height–width ratio of your stitch. The general formula of adding $2\pi H/W$ stitches allows you to check the correct number of stitches needed for any type of stitch.

But what if you wanted to crochet a sphere? Unlike the plane, or any flat surface with zero curvature, a sphere has constant positive curvature and curves back onto itself instead of stretching to infinity. To crochet an approximation of a sphere, we need spherical geometry. 

\begin{figure}[t]
\centering
\begin{minipage}[b]{0.95\textwidth} 
	\includegraphics[width=\textwidth]{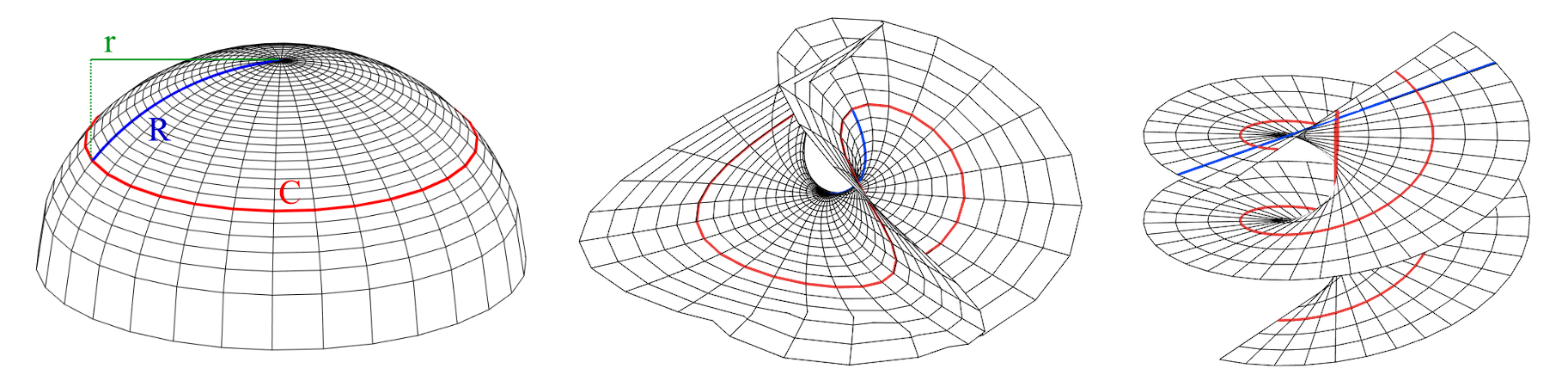}
\end{minipage}
\caption{To crochet a sphere (left), or an Enneper surface (middle) we need to consider the intrinsic radius $R$ and the circumference of a circle $C$. For the helicoid (right) $R=r$ and instead of circumference one needs to consider the length of rows.}
\label{Fig:MatlabFigures}
\end{figure}

Consider a sphere with radius $S$, and draw a circle with intrinsic radius $R$ (the radius measured along the surface of the sphere), as shown in Figure \ref{Fig:MatlabFigures}. The circumference of this circle is
\begin{equation*}
C = 2\pi S \cdot \sin\bigg(\frac{R}{S}\bigg).
\end{equation*}
Unlike in the case of a disc, the circumference does not grow linearly. Notice that $S/R \cdot \sin(R/S) \leq 1$, which, rather naturally, implies that to create a sphere you should add fewer stitches per round than when creating a flat circle. Also, the largest possible circle, with intrinsic radius $R = \pi S$, divides the sphere into two equal halves. This means that to crochet a sphere with radius $S$, we only need to calculate the required stitches up to the level $\ell \leq \pi S / H$, after which the second half of the sphere can be completed as the first part but in reverse. You can now crochet a close approximation of a sphere of any size by using the formula
\begin{align*}
\frac{2\pi S}{W} \cdot \sin\bigg(\frac{H \cdot \ell}{S}\bigg)
\end{align*}
to calculate the number of required stitches on round $\ell$.

Discs and spheres have been crocheted for a long time, but the idea of crocheting a hyperbolic plane, which has constant negative curvature, was not introduced until 1997 by Daina Taimina \cite{Henderson2001}. One way of constructing a hyperbolic surface is to consider the circumference of a circle drawn on a hyperbolic plane:
\begin{equation*}
C = 2\pi S \cdot \sinh\bigg(\frac{R}{S}\bigg),
\end{equation*}
where $S$ is the radius of the hyperbolic plane. Note that $S/R \cdot \sinh(R/S) \geq 1$, which means that the hyperbolic plane requires more additional stitches per round than the flat disc. Detailed instructions for crocheting hyperbolic planes can be found in \cite{Taimina2018}.

\subsection{Minimal surfaces} 

The plane, sphere, and hyperbolic plane all have constant Gaussian curvature, meaning that they look the same at every point, allowing general formulas for the circumference of a circle. This makes calculating crochet models simpler, but since we crochet in rounds (or rows) we do not need to know the circumference of every possible circle, only the circumferences of the circles centred on the starting point. We will next consider a group of more general ‘crochet-symmetric’ shapes that allow similarly simple crochet instructions (at least up to intersections), requiring only one type of stitch, with the additional or removed stitches spaced evenly around a round.

Imagine dipping a wire frame into soapy water. The soap film will form an optimal surface that minimises the surface area bounded by the wire. In mathematics, such area-minimising ‘soap film’ surfaces are called minimal surfaces. More generally, we define minimal surfaces to be surfaces that locally minimise their area. These surfaces do not necessarily minimise the surface area globally, they are allowed to self-intersect, and they do not have to have a boundary. Minimal surfaces have the interesting property that their mean curvature, given by the sum of the principal curvatures, is zero. This means that the principal curvatures must have the same magnitude but opposite signs at every point, so the Gaussian curvature is either negative or zero. Unlike the plane, sphere, and hyperbolic plane, minimal surfaces do not have constant Gaussian curvature.

The study of minimal surfaces began in 1761, when the Italian mathematician Joseph-Louis Lagrange addressed the task of finding a surface with the smallest possible area given a closed boundary \cite{Lagrange1761}. Lagrange proposed that such minimal surfaces could locally be represented as graphs of solutions to the so-called minimal surface equation, but the only concrete example he could provide was the plane.
In 1873, Joseph Plateau published his experimental findings on minimal surfaces formed by soap films \cite{Plateau1873}. His work offered a clear physical interpretation of the problem and helped spread interest in it beyond mathematics. Consequently, the problem of demonstrating the existence of a minimal surface for a given boundary became known as Plateau’s problem. Only special cases of the problem were solved before the 1930s, when the general case was proved independently by Jesse Douglas \cite{Douglas1931} and Tibor Radó \cite{Rado1930}.

%\begin{wrapfigure}{r}{6cm}
%\includegraphics[width=6cm]{figures/HelicoidCatenoid}
%\caption{Crocheted catenoid (top) and helicoid (bottom).}
%\label{Fig:HelicoidCatenoid}
%\end{wrapfigure}

\begin{figure}[!b]
\centering
\begin{minipage}[b]{0.95\textwidth} 
	\includegraphics[width=\textwidth]{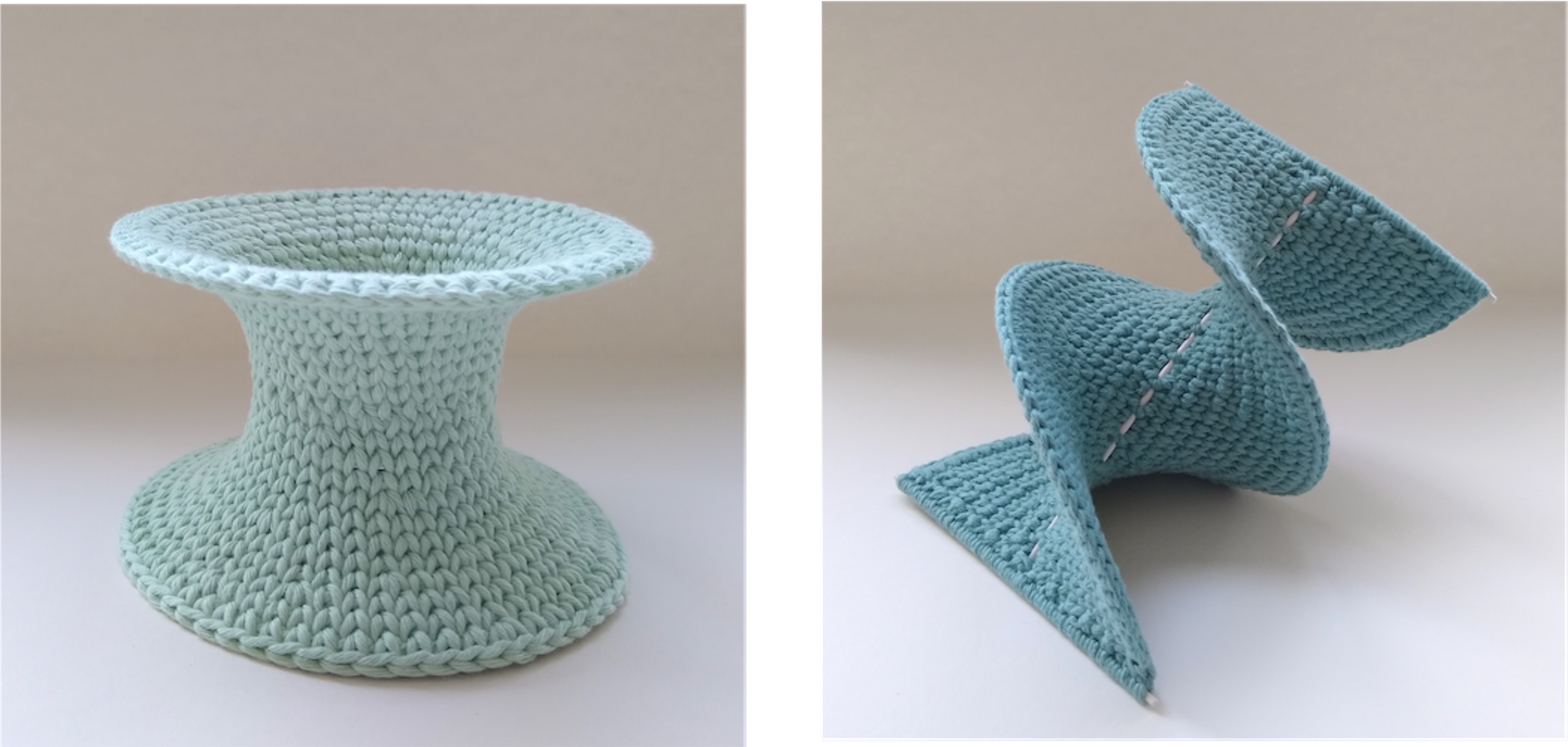}
\end{minipage}
\caption{Crocheted catenoid (left) and helicoid (right).}
\label{Fig:HelicoidCatenoid}
\end{figure} 

Before the 19\textsuperscript{th} century, the plane, catenoid, and helicoid were the only known examples of minimal surfaces. The catenoid was first described by the Swiss mathematician Leonhard Euler in 1744. A catenoid can be created by rotating a catenary curve, a curve formed by a heavy freely hanging rope with fixed ends, around an axis. Note that a catenary curve is not the same as a parabola, even though the two look quite similar. The Gaussian curvature of a point on the catenoid depends only on its distance from the origin. The catenoid is strongly negatively curved near its neck, but the further you move from the centre, the flatter the surface becomes. Because of its rotational symmetry, you can crochet a catenoid following a simple pattern that uses only one type of stitch, with the increases evenly spaced across each round, as shown in Figure \ref{Fig:HelicoidCatenoid}. Unfortunately, the catenoid is the only minimal surface that can be created by rotating a curve.

The next minimal surface, the right helicoid, was described by the French mathematician and engineer Jean-Baptiste Meusnier in 1776. Meusnier also showed that the catenoid satisfies Lagrange’s condition and gave a geometric interpretation of Lagrange’s minimal surface equation, noting that it is equivalent to the surface having zero mean curvature \cite{Meusnier1785}.
The helicoid can be constructed by connecting the two branches of a double helix with horizontal lines passing through the vertical centre line, as shown in Figure \ref{Fig:MatlabFigures}. We can draw such a straight line through any point on the helicoid, meaning that it is a ruled surface. In fact, the helicoid and the plane are the only ruled minimal surfaces.

Even though the catenoid and helicoid look quite different, they are isometric. This means they share the same intrinsic geometry, so a helicoid can be twisted into a catenoid (that is, cut along one catenary). It also means that you only need to calculate the crochet pattern for one of them. See \cite{Taimina2018} for crochet instructions for both the catenoid and helicoid.

As mentioned above, the catenoid is the only minimal surface that is also a surface of revolution. However, there exists a group of minimal surfaces that are intrinsically surfaces of revolution, which is sufficient to allow simple crochet instructions. These minimal surfaces were introduced in 1862 by the French mathematician and engineer Edmond Bour \cite{Bour1862}, in a paper published a few years before the famous Enneper–Weierstrass representation, which provided the first general method for finding minimal surfaces. The catenoid and helicoid cannot be expressed in the form used by Bour, but they can be included in this group when considered via their Enneper–Weierstrass data.

In what follows, we examine three different families of surfaces: Enneper’s, Richmond’s, and Bour’s surfaces. All three can be derived from the set of equations Bour presented in his 1862 paper, but Enneper’s surfaces are named after Alfred Enneper, who described them in greater detail in 1864, and Richmond’s surfaces after Herbert William Richmond, who studied them in 1904. The final example we discuss is usually referred to simply as Bour’s minimal surface. For more on the history and recent advances in the study of minimal surfaces, see \cite{Perez2017}, and for a comprehensive introduction, see \cite{Nitsche1989}.

\begin{figure}[t]
\centering
\begin{minipage}[b]{0.95\textwidth} 
	\includegraphics[width=\textwidth]{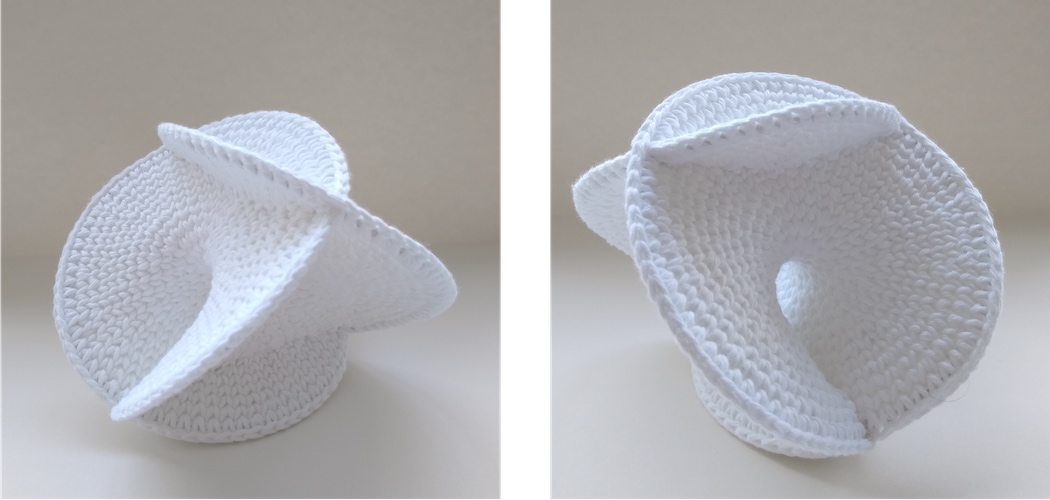}
\end{minipage}
\caption{Two views of the self intersectiong Enneper's surface.}\label{Fig:Enneper1}
\end{figure}

We start with Enneper's surfaces, illustrated in Figure \ref{Fig:Enneper1}, which are given by the equations
\begin{align*}\label{eg:Enneper}
\begin{split}
x(r,\theta) &= r\cos(\theta) - \frac{r^{2n-1}}{2n-1}\cos\big((2n-1)\theta\big)\\
y(r,\theta) &= r\sin(\theta) + \frac{r^{2n-1}}{2n-1}\sin\big((2n-1)\theta\big)\\
z(r,\theta) &= \frac{2r^n}{n}\cos(n\theta),
\end{split}
\end{align*}
where $r \geq 0$ gives the distance from the origin, $0 \leq \theta \leq 2\pi$ the rotation angle, and $n = 2, 3, \dots$ the rotational symmetry of the surface. Choosing $n = 1$ produces a flat disc, which is a trivial minimal surface. The model can of course be scaled to any desired size.

To create an accurate crocheted Enneper’s surface, we need to calculate the circumference of a circle centred at the origin with an intrinsic radius $R = \ell \cdot H$. This requires solving arc lengths of curves, which can be done in two steps:
\begin{itemize}
\item[1.] Calculate the circumference $C$ of a circle as a function of $r$.
\item[2.] Determine $r$ as a function of the intrinsic radius $R$.
\end{itemize}
In the first step, we determine the arc length of the red curve in Figure \ref{Fig:MatlabFigures} for a fixed $r$. Note that the $r$ in the above equations is not the intrinsic radius $R$, which is measured along the surface. However, when adding a new round to the model, the intrinsic radius $R$ grows by the height of the stitch, which means we also need to calculate the arc length of the blue curve in Figure \ref{Fig:MatlabFigures}. For Enneper’s surfaces we obtain
\begin{align*}
C(r) & = 2\pi(r+r^{2n-1})\\
R(r) & = r+\frac{r^{2n-1}}{2n-1}.
\end{align*}
Additionally, the length of a line segment on a circle $C$ between angles $\theta_1$ and $\theta_2$ depends only on the difference between $\theta_1$ and $\theta_2$. This allows the increases to be spaced evenly on every round. The intrinsic radius $R$ does not depend on the angle $\theta$, so the model can be made using only one type of stitch with height $H$.

To find the number of required stitches on round $\ell$, we first solve $r$ as a function of $R = \ell \cdot H$. In general, this must be done numerically. We then calculate the circumference of the circle with intrinsic radius $\ell \cdot H$ and divide this by the width of a stitch $W$. That is, the number of stitches on round $\ell$ is $$2\pi(r(\ell)+r(\ell)^{2n-1})/W$$ rounded to the nearest integer.

\begin{figure}[t]
\centering
\begin{minipage}[b]{0.95\textwidth} 
\includegraphics[width=\textwidth]{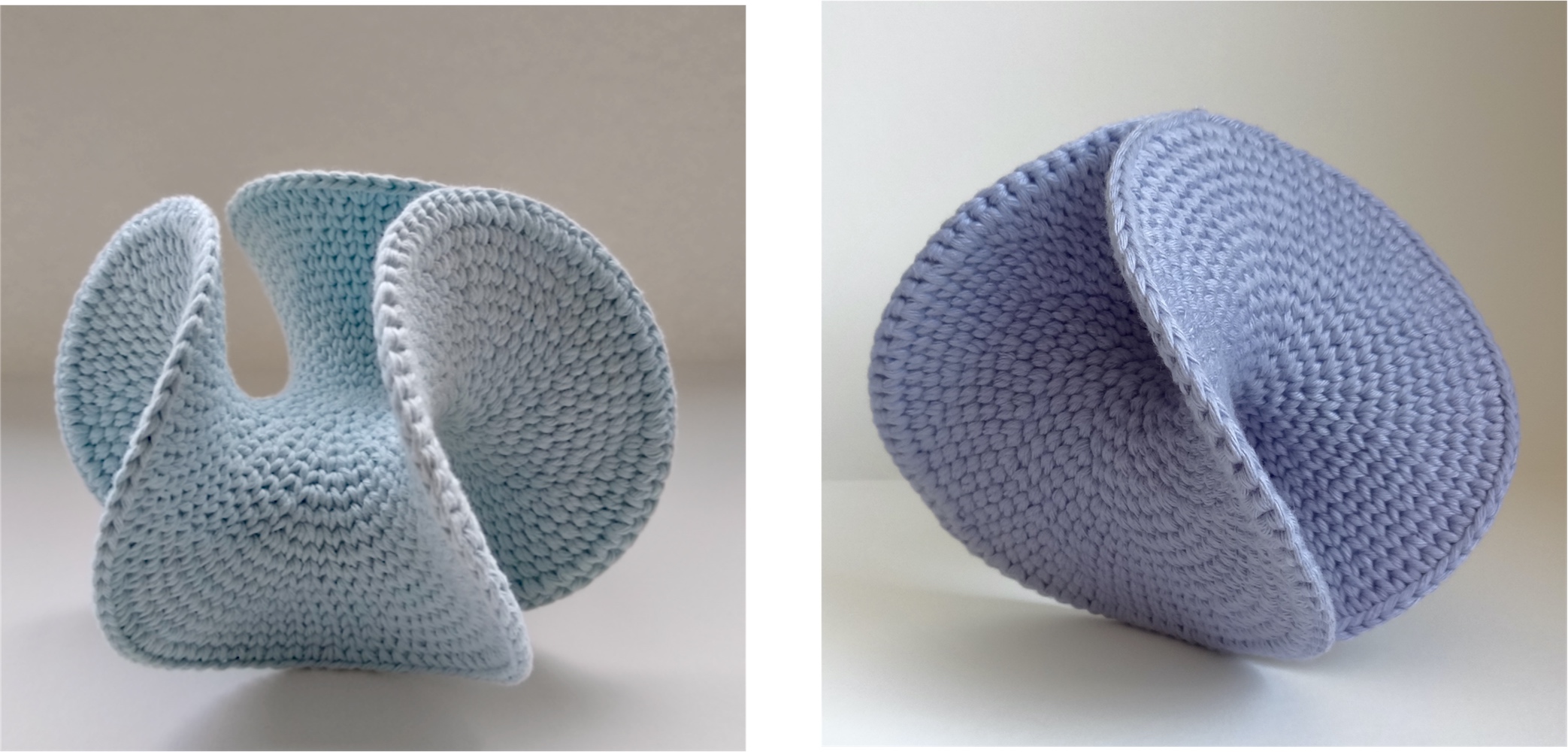}
\end{minipage}
\caption{On the left, a crocheted Enneper surface with 3-fold symmetry; on the right an Enneper’s surface with 2-fold symmetry crocheted to just before the intersection.}
\label{Fig:Enneper2}
\end{figure}

Close to the origin, a small piece of an Enneper’s surface with order-2 symmetry resembles a saddle, but as $r$ increases the surface begins to intersect itself. The intersection does not occur along a straight line, and some stitches must be shifted from one side to the other, adding a degree of challenge to the model. Table \ref{Tab:Enneper} provides instructions for crocheting an order-2 Enneper’s surface up to the point of intersection, as shown on the right in Figure \ref{Fig:Enneper2}. For a more detailed explanation of the mathematics behind the models, as well as instructions for a version with intersections, see \cite{Kekkonen2023}.

When crocheting mathematical models, I recommend using a yarn with minimal stretch and a crochet hook smaller than the manufacturer’s recommendation to achieve a tighter gauge. For the minimal surfaces shown in this book, I used 10-ply cotton yarn (50 g $\approx$ 75 m) with a 3 mm crochet hook. Begin by making a test tube to measure the height and width of your stitches; because the model is worked in rounds, a tube will give more accurate measurements than a traditional flat swatch. The models are quite sensitive to stitch size, so aim for consistency. You can find instructions for three different stitch heights in Table \ref{Tab:Enneper}. Each model requires just over 1 500 stitches and produces a piece about 13 cm wide, with the top and bottom just meeting. Using thicker (or thinner) yarn while keeping the height–width ratio constant will yield a larger (or smaller) model. You can use either the traditional single crochet stitch or the split stitch for a smoother, more rigid surface. If you are new to crochet, I suggest starting with the traditional single stitch, as it is easier to learn and gentler on the wrists.

Below I describe how to crochet the Enneper’s surface in the case where the width and height of your crochet stitch are the same. To begin, create a magic loop and crochet six stitches into it, then close the round with a slip stitch. On the next round ($\ell = 2$), add eight stitches evenly spaced and close the round with a slip stitch. At this point, you should have a total of 14 stitches. Continue working on the model, always distributing the added stitches evenly on each round. Whenever possible, try to stagger your increases so that they are not stacked on top of each other.You can also crochet continuously without closing each round with a slip stitch, but be sure to mark the start of the round. The final round may have a slightly less neat finish when using continuous rounds, but this technique helps avoid seams that might appear when closing the rounds with a slip stitch. I recommend crocheting a thin plastic cable (e.g., transparent trimmer line) into your model during the last round to provide support and help maintain the correct shape. After completing the final round, divide the stitches into four equal parts (50 stitches each) and sew the top and bottom together. Weave in the ends and your model is finished. If you prefer a smaller Enneper surface, you can simply stop at any round.

\newpage

As previously mentioned, the model is quite sensitive to an incorrect gauge. If you find that the edge of your model does not curve enough to meet at the last round, the cause is likely that your stitches are too tall, assuming the width remains constant. This may seem counterintuitive, but taller stitches increase the radius of the circle while the circumference remains the same, which results in a surface that is less curved. Conversely, if it feels like there is not enough space for the final round, your stitches are probably too short, creating a model with too much curvature.

\begin{table}
\begin{center}
\begin{tabular}{|c|c|c|}
\cline{1-3}
\multicolumn{3}{|c|}{ H $=0.4$ cm, W $=0.5$ cm}\\
\hline $\ell$ & \hspace*{1.2mm} $\Delta$ N \hspace*{1.2mm} & \hspace*{2.5mm} N \hspace*{2.5mm} \\
\hline 1 & & 5 \\
\hline 2 & 6 & 11 \\
\hline 3 & 7 & 18 \\
\hline 4 & 7 & 25 \\
\hline 5 & 9 & 34 \\
\hline 6 & 9 & 43 \\
\hline 7 & 10 & 53 \\
\hline 8 & 10 & 63 \\
\hline 9 & 11 & 74 \\
\hline 10 & 11 & 85 \\
\hline 11 & 11 & 96 \\
\hline 12 & 11 & 107 \\
\hline 13 & 12 & 119 \\
\hline 14 & 12 & 131 \\
\hline 15 & 12 & 143 \\
\hline 16 & 12 & 155 \\
\hline 17 & 12 & 167 \\
\hline 18 & 13 & 180 \\
\hline  &  &  \\
\hline \text {Stitches} & & 1512
\\
\hline
\end{tabular}
\quad
\begin{tabular}{|c|c|c|}
\cline{1-3}
\multicolumn{3}{|c|}{H $=0.45$ cm, W $=0.5$ cm}\\
\hline $\ell$ & \hspace*{1.2mm} $\Delta$ N \hspace*{1.2mm} & \hspace*{2.5mm} N \hspace*{2.5mm}\\
\hline 1 & & 6 \\
\hline 2 & 6 & 12 \\
\hline 3 & 8 & 20 \\
\hline 4 & 9 & 29 \\
\hline 5 & 10 & 39 \\
\hline 6 & 10 & 49 \\
\hline 7 & 12 & 61 \\
\hline 8 & 12 & 73 \\
\hline 9 & 12 & 85 \\
\hline 10 & 12 & 97 \\
\hline 11 & 13 & 110 \\
\hline 12 & 13 & 123 \\
\hline 13 & 13 & 136 \\
\hline 14 & 14 & 150 \\
\hline 15 & 14 & 164 \\
\hline 16 & 14 & 178 \\
\hline 17 & 14 & 192 \\
\hline  &  &  \\
\hline  &  &  \\
\hline \text {Stitches} & & 1525
\\
\hline
\end{tabular}
\quad
\begin{tabular}{|c|c|c|}
\cline{1-3}
\multicolumn{3}{|c|}{H $=0.5$ cm, W $=0.5$ cm}\\
\hline $\ell$ & \hspace*{1.2mm} $\Delta$ N \hspace*{1.2mm} & \hspace*{2.5mm} N \hspace*{2.5mm} \\
\hline 1 & & 6 \\
\hline 2 & 8 & 14 \\
\hline 3 & 9 & 23 \\
\hline 4 & 10 & 33 \\
\hline 5 & 11 & 44 \\
\hline 6 & 12 & 56 \\
\hline 7 & 13 & 69 \\
\hline 8 & 13 & 82 \\
\hline 9 & 14 & 96 \\
\hline 10 & 14 & 110 \\
\hline 11 & 14 & 124 \\
\hline 12 & 15 & 139 \\
\hline 13 & 15 & 154 \\
\hline 14 & 15 & 169 \\
\hline 15 & 15 & 184 \\
\hline 16 & 16 & 200 \\
\hline  &  &  \\
\hline  &  &  \\
\hline  &  &  \\
\hline \text {Stitches} & & 1504
\\
\hline
\end{tabular}
\caption{The number added stitches $\Delta N$ and number of stitches N at round $\ell$ for models with gauge W $=0.5$ cm and H $=0.4$ cm, H $=0.45$ cm or H $=0.5$ cm. All the models require a little bit over $1500$ stitches to complete.}
\label{Tab:Enneper} 
\end{center}
\end{table}

\vspace{-2mm}
\begin{wrapfigure}{r}{5.5cm}
\includegraphics[width=5.5cm]{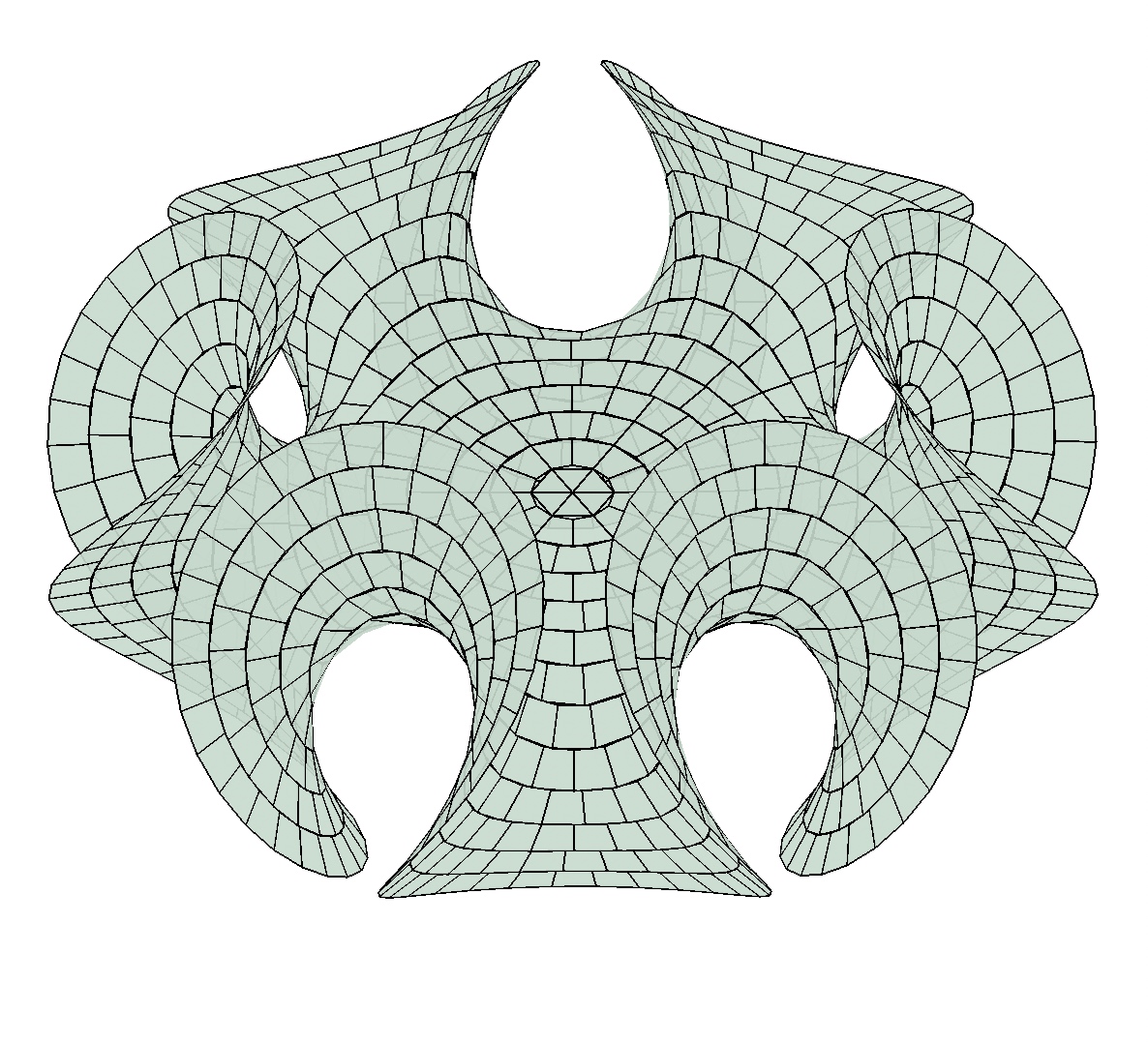}
\caption{Enneper’s surface with 6-fold symmetry, showing the distinct almost flat disc at the centre.}
\label{Fig:Enneper6Fold}
\end{wrapfigure}

We can also consider Enneper’s surfaces with higher-order rotational symmetries. For example, Figure \ref{Fig:Enneper2} shows a crocheted Enneper surface with three wavy ends instead of two. In the order-two case, the surface has a large negative Gaussian curvature at the origin, with the curvature tending towards zero as we move further away from the centre. For higher-order symmetries, however, the curvature at the origin is zero, and the larger the order, the more pronounced the appearance of an almost flat central disc becomes, as can be seen in Figure \ref{Fig:Enneper6Fold}. This central region is surrounded by waves where the Gaussian curvature reaches its peak before gradually flattening out again.

\clearpage

\begin{figure}[H]
\centering
\begin{minipage}[b]{0.95\textwidth} 
\includegraphics[width=\textwidth]{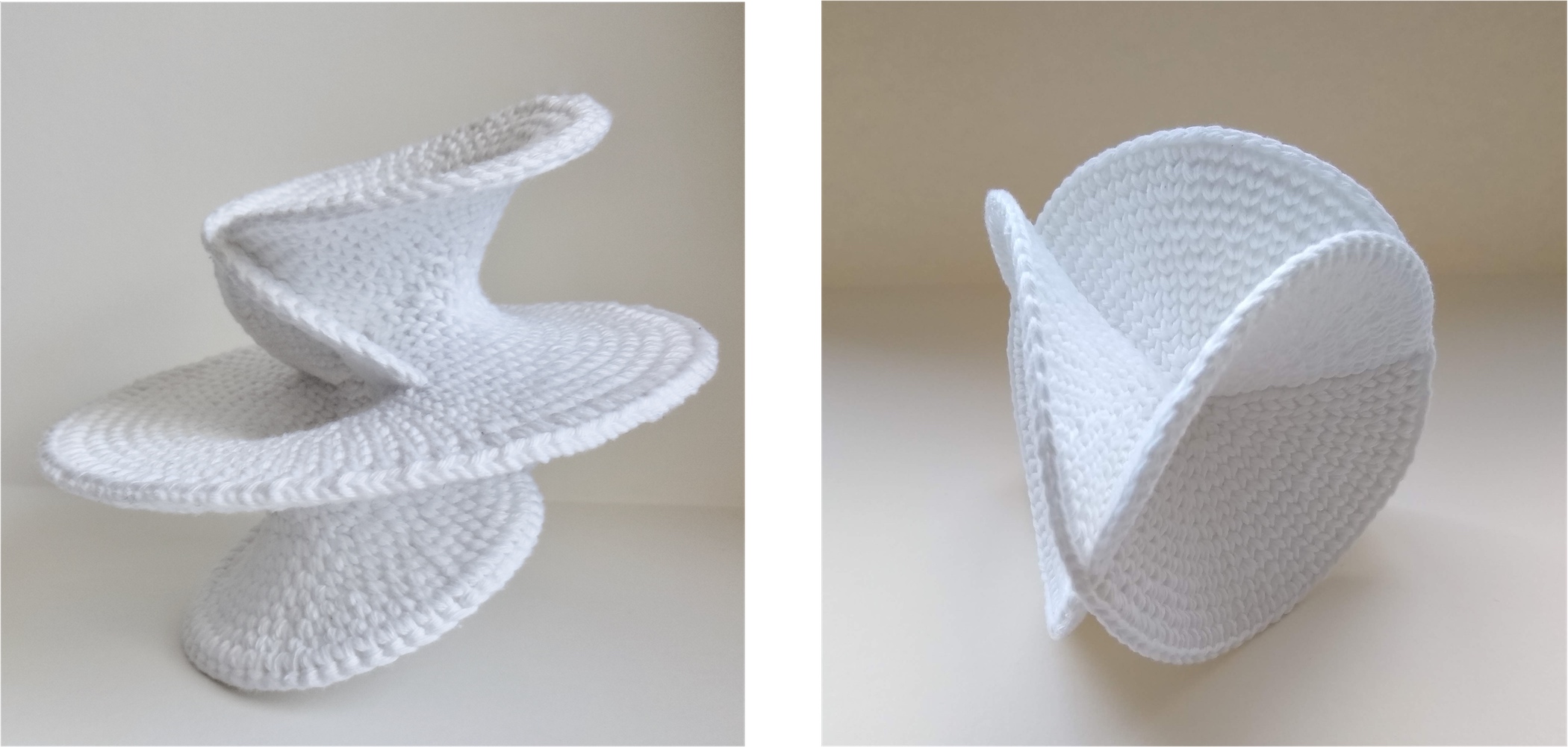}
\end{minipage}
\caption{Crocheted Richmond's minimal surface (left) and Bour's minimal surface (right).}
\label{Fig:BourRichmond}
\end{figure}

Another example of minimal surfaces that are intrinsically surfaces of revolution is the family of Richmond's surfaces, parametrised by
\begin{align*}
x(r,\theta) & = -\frac{\cos(\theta)}{r}-\frac{r^{2n+1}}{2n+1}\cos((2n+1)\theta)\\
y(r,\theta) & = -\frac{\sin(\theta)}{r}-\frac{r^{2n+1}}{2n+1}\sin((2n+1)\theta)\\
z(r,\theta) & = \frac{2r^n}{n}\cos(n\theta),
\end{align*}
where $r>0$, $0\leq\theta\leq2\pi$, and $n=1,2,\dots$. Richmond's surfaces have one Enneper-type end at the middle and one planar-type end stretching outward, which is why they are sometimes referred to as planar Enneper surfaces. In practice, we consider models where $r \in [r_1, r_2]$. Choosing a smaller $r_1$ allows the planar edge to extend further, whereas a larger $r_2$ produces a more prominent intersection area in the middle.

Since a Richmond’s surface extends in two directions, there is no single best place to begin the crochet model. Because removing stitches usually looks less tidy than adding them, one option is to start from the shortest round and work in both directions separately. The intersections of Richmond’s surfaces are also more complex than those of Enneper’s surfaces. From a crocheting perspective, there are two distinct types of self-intersections: the Enneper-type, where a round crosses itself several times, and a straight-line intersection, where the surface appears to pass directly through an already completed section. This second type can be achieved by beginning the model from the first round before the intersection and crocheting in both directions alternatingly, allowing the two parts to intersect naturally as was done to achieve the model shown in Figure \ref{Fig:BourRichmond}.

The last surface we consider is Bour’s minimal surface, which can be written as
\begin{align*}
x(r,\theta) & = r\cos(\theta)-\frac{1}{2}r^2\cos(2\theta)\\
y(r,\theta) & = -r\sin(\theta)-\frac{1}{2}r^2\sin(2\theta)\\
z(r,\theta) & = \frac{4}{3}r^{\frac{3}{2}}\cos\bigg(\frac{3}{2}\theta\bigg),
\end{align*}
where $r\in[0,1]$ and $\theta\in[0,4\pi]$.
The surface intersects itself along three rays lying on a plane that divides it into two symmetric halves. These rays meet at equal angles at the origin, partitioning the surface into six identical sections. Crocheting this surface can be tricky to begin with, as it requires a dual magic ring construction, but once the initial, less stable rounds are complete, the work becomes straightforward. Unlike the Enneper and Richmond surfaces, the six sections are separated by straight lines, so stitches remain within their section and do not need to be shifted from one part to another.

%%%%%%%%%%%%%%%%%%%%%%%%%%%%%%%%%%%%%%
% References %
\bibliographystyle{aomplain}
\bibliography{ReferencesRecreational}

\begin{thebibliography}{34}
\providecommand{\natexlab}[1]{#1}
\providecommand{\url}[1]{\texttt{#1}}
\expandafter\ifx\csname urlstyle\endcsname\relax
  \providecommand{\doi}[1]{doi: #1}\else
  \providecommand{\doi}{doi: \begingroup \urlstyle{rm}\Url}\fi

\bibitem[Adams(2004)]{Adams2004}
Colin~C. Adams.
\newblock \emph{The knot book}.
\newblock American Mathematical Society, Providence, RI, 2004.
\newblock ISBN 0-8218-3678-1.
\newblock URL
  \url{https://www.math.cuhk.edu.hk/course_builder/1920/math4900e/Adams--The%20Knot%20Book.pdf}.

\bibitem[al-Razz\~az al-Jazar\~i (translated and annotated~by
  D.~R.~Hill)(1974)]{Ibn1206}
Ibn al-Razz\~az al-Jazar\~i (translated and annotated~by D.~R.~Hill).
\newblock \emph{The Book of Knowledge of Ingenious Mechanical Devices}.
\newblock Springer Dordrecht, 1974.
\newblock ISBN 9027-70-329-9.
\newblock \doi{10.1007/978-94-010-2573-7}.
\newblock URL \url{https://doi.org/10.1007/978-94-010-2573-7}.

\bibitem[Ashton et~al.(2011{\natexlab{a}})Ashton, Cantarella, Piatek, and
  Rawdon]{Ashton11b}
Ted Ashton, Jason Cantarella, Michael Piatek, and Eric~J. Rawdon.
\newblock Atlas of tight links, 2011{\natexlab{a}}.
\newblock URL
  \url{https://jasoncantarella.com/downloads/papers/rrpaper/TightKnotCatalogue.pdf}.

\bibitem[Ashton et~al.(2011{\natexlab{b}})Ashton, Cantarella, Piatek, and
  Rawdon]{Ashton2011}
Ted Ashton, Jason Cantarella, Michael Piatek, and Eric~J. Rawdon.
\newblock Knot tightening by constrained gradient descent.
\newblock \emph{Exp. Math.}, 20\penalty0 (1):\penalty0 57--90,
  2011{\natexlab{b}}.
\newblock ISSN 1058-6458.
\newblock \doi{10.1080/10586458.2011.544581}.
\newblock URL \url{https://doi.org/10.1080/10586458.2011.544581}.

\bibitem[Bour(1862)]{Bour1862}
Edmond Bour.
\newblock Theorie de la deformation des surfaces.
\newblock \emph{Journal de l'Ecole Polytechnique}, 22(39):\penalty0 1--148,
  1862.
\newblock URL \url{https://gallica.bnf.fr/ark:/12148/bpt6k433694t/f5.item}.

\bibitem[Buck and Flapan(2009)]{Buck2009}
Dorothy Buck and Erica Flapan.
\newblock \emph{Applications of Knot Theory: American Mathematical Society,
  Short Course, January 4-5, 2008, San Diego, California}, volume~66.
\newblock American Mathematical Soc., 2009.
\newblock ISBN 978-0-8218-4466-3.

\bibitem[Cartwright and Gonz{\'a}lez(2016)]{Cartwright2016}
Julyan H.~E. Cartwright and Diego~L. Gonz{\'a}lez.
\newblock M{\"o}bius strips before m{\"o}bius: Topological hints in ancient
  representations.
\newblock \emph{The Mathematical Intelligencer}, 38\penalty0 (2):\penalty0
  69--76, 2016.
\newblock \doi{10.1007/s00283-016-9631-8}.
\newblock URL \url{https://doi.org/10.1007/s00283-016-9631-8}.

\bibitem[Crum~Brown(1886)]{Brown1886}
Alexander Crum~Brown.
\newblock On a case of interlacing surfaces.
\newblock \emph{Proceedings of the Royal Society of Edinburgh}, 13:\penalty0
  382--386, 1886.

\bibitem[Denne et~al.(2006)Denne, Diao, and Sullivan]{Denne2006}
Elizabeth Denne, Yuanan Diao, and John~M. Sullivan.
\newblock Quadrisecants give new lower bounds for the ropelength of a knot.
\newblock \emph{Geom. Topol.}, 10:\penalty0 1--26, 2006.
\newblock ISSN 1465-3060.
\newblock \doi{10.2140/gt.2006.10.1}.
\newblock URL \url{https://doi.org/10.2140/gt.2006.10.1}.

\bibitem[Dong(2023)]{Dong2023}
Shiying Dong.
\newblock Sculpting mapping cylinders: Seamless crochet of topological
  surfaces.
\newblock In \emph{Proceedings of Bridges 2023: Mathematics, Art, Music,
  Architecture, Culture}, pages 559--566. Tessellations Publishing, 2023.
\newblock ISBN 978-1-938664-45-8.
\newblock URL
  \url{http://archive.bridgesmathart.org/2023/bridges2023-559.html}.

\bibitem[Douglas(1931)]{Douglas1931}
Jesse Douglas.
\newblock Solution of the problem of {P}lateau.
\newblock \emph{Transactions of the American Mathematical Society}, 33\penalty0
  (1):\penalty0 263--321, 1931.
\newblock URL \url{https://doi.org/10.2307/1989472}.

\bibitem[Dunning(2015)]{Dunning2015}
David~E. Dunning.
\newblock What are models for? {A}lexander {C}rum {B}rown's knitted
  mathematical surfaces.
\newblock \emph{Math. Intelligencer}, 37\penalty0 (2):\penalty0 62--70, 2015.
\newblock ISSN 0343-6993.
\newblock \doi{10.1007/s00283-014-9480-2}.
\newblock URL \url{https://doi.org/10.1007/s00283-014-9480-2}.

\bibitem[Henderson and Taimi\c{n}a(2001)]{Henderson2001}
David~W. Henderson and Daina Taimi\c{n}a.
\newblock Crocheting the hyperbolic plane.
\newblock \emph{Math. Intelligencer}, 23\penalty0 (2):\penalty0 17--28, 2001.
\newblock ISSN 0343-6993.
\newblock URL \url{https://doi.org/10.1007/BF03026623}.

\bibitem[Horner et~al.(2016)Horner, Miller, Steed, and Sutcliffe]{Horner2016}
Kate~E Horner, Mark~A Miller, Jonathan~W Steed, and Paul~M Sutcliffe.
\newblock Knot theory in modern chemistry.
\newblock \emph{Chemical Society Reviews}, 45\penalty0 (23):\penalty0
  6432--6448, 2016.

\bibitem[Inman(1903)]{Inman1903}
E.~R. Inman.
\newblock Half twisting a quarter twist belt.
\newblock \emph{American Machinist}, 26:\penalty0 4466--4467, 1903.
\newblock URL
  \url{https://archive.org/details/sim_american-machinist_june-04-december-31-1903_26_23-52/page/1466/}.

\bibitem[Irving(2004)]{Irving2004}
Claire Irving.
\newblock Making the real projective plane.
\newblock \emph{The Mathematical Gazette}, 89\penalty0 (516):\penalty0
  417--423, 2004.
\newblock URL \url{https://www.jstor.org/stable/3621933}.

\bibitem[Kekkonen(2023)]{Kekkonen2023}
Hanne Kekkonen.
\newblock Crocheting {B}our's $\mathcal{B}_m$ minimal surfaces.
\newblock \emph{arXiv preprint arXiv:2306.15378}, 2023.
\newblock URL \url{https://arxiv.org/pdf/2306.15378.pdf}.

\bibitem[Lagrange(1761)]{Lagrange1761}
Joseph-{L}ouis Lagrange.
\newblock Essai d'une nouvelle methode pour de'terminer les maxima et les
  minima des formules integrales indefinies.
\newblock \emph{Miscellanea Taurinensia}, t. II:\penalty0 335--362, 1761.
\newblock URL \url{https://gallica.bnf.fr/ark:/12148/bpt6k2155691/f385}.

\bibitem[Luotoniemi(2019)]{Luotoniemi2019}
Taneli Luotoniemi.
\newblock \emph{{Hyperspatial Interlace - Grasping Four-dimensional Geometry
  through Crafted Models}}.
\newblock Doctoral thesis, School of Arts, Design and Architecture, 2019.
\newblock URL \url{http://urn.fi/URN:ISBN:978-952-60-8480-0}.

\bibitem[Meusnier(1785 (presented 1776))]{Meusnier1785}
Jean~Baptiste Meusnier.
\newblock M{\'e}moire sur la courbure des surfaces.
\newblock \emph{Mem des savan etrangers}, 10:\penalty0 477--510, 1785
  (presented 1776).

\bibitem[Nitsche(1989)]{Nitsche1989}
Johannes~C Nitsche.
\newblock \emph{Lectures on minimal surfaces: vol. 1}.
\newblock Cambridge university press, 1989.

\bibitem[Osinga and Krauskopf(2004)]{Osinga2004}
Hinke~M. Osinga and Bernd Krauskopf.
\newblock Crocheting the {L}orenz manifold.
\newblock \emph{Math. Intelligencer}, 26\penalty0 (4):\penalty0 25--37, 2004.
\newblock ISSN 0343-6993.
\newblock \doi{10.1007/BF02985416}.
\newblock URL \url{https://doi.org/10.1007/BF02985416}.

\bibitem[P\'{e}rez(2017)]{Perez2017}
Joaqu\'{\i}n P\'{e}rez.
\newblock A new golden age of minimal surfaces.
\newblock \emph{Notices Amer. Math. Soc.}, 64\penalty0 (4):\penalty0 347--358,
  2017.
\newblock ISSN 0002-9920.
\newblock \doi{10.1090/noti1500}.
\newblock URL \url{https://doi.org/10.1090/noti1500}.

\bibitem[Plateau(1873)]{Plateau1873}
Joseph Antoine~Ferdinand Plateau.
\newblock \emph{Statique exp{\'e}rimentale et th{\'e}orique des liquides soumis
  aux seules forces mol{\'e}culaires}, volume~2.
\newblock Gauthier-Villars, 1873.

\bibitem[Rad{\'o}(1930)]{Rado1930}
Tibor Rad{\'o}.
\newblock On {P}lateau's problem.
\newblock \emph{Annals of Mathematics}, pages 457--469, 1930.
\newblock URL \url{https://doi.org/10.2307/1968237}.

\bibitem[Rawdon(2003)]{Rawdon2003}
Eric~J. Rawdon.
\newblock Can computers discover ideal knots?
\newblock \emph{Experiment. Math.}, 12\penalty0 (3):\penalty0 287--302, 2003.
\newblock ISSN 1058-6458.
\newblock URL \url{http://projecteuclid.org/euclid.em/1087329232}.

\bibitem[Reid(1971)]{Reid1971}
Miles Reid.
\newblock The knitting of surfaces.
\newblock \emph{Eureka - The Journal of the {A}rchimedeans}, 34:\penalty0
  21--26, 1971.
\newblock URL \url{http://homepages.warwick.ac.uk/~masda/knit_surfaces.pdf}.

\bibitem[Seifert(1935)]{Seifert1935}
H.~Seifert.
\newblock \"{U}ber das {G}eschlecht von {K}noten.
\newblock \emph{Math. Ann.}, 110\penalty0 (1):\penalty0 571--592, 1935.
\newblock ISSN 0025-5831.
\newblock \doi{10.1007/BF01448044}.
\newblock URL \url{https://doi.org/10.1007/BF01448044}.

\bibitem[Singal et~al.(2023)Singal, Dimitriyev, Gonzalez, Quinn, and
  Matsumoto]{Singal2023}
Krishma Singal, Michael~S Dimitriyev, Sarah~E Gonzalez, Sam Quinn, and
  Elisabetta~A Matsumoto.
\newblock Programming mechanics in knitted materials, stitch by stitch.
\newblock \emph{arXiv preprint arXiv:2302.13467}, 2023.
\newblock URL \url{https://arxiv.org/pdf/2302.13467.pdf}.

\bibitem[Stasiak et~al.(1998)Stasiak, Katritch, and Kauffman]{Stasiak1998}
A.~Stasiak, V.~Katritch, and L.~H. Kauffman, editors.
\newblock \emph{Ideal knots}, volume~19 of \emph{Series on Knots and
  Everything}.
\newblock World Scientific Publishing Co., Inc., River Edge, NJ, 1998.
\newblock ISBN 981-02-3530-5.
\newblock \doi{10.1142/9789812796073}.
\newblock URL \url{https://doi.org/10.1142/9789812796073}.

\bibitem[Taimi\c{n}a(2018)]{Taimina2018}
Daina Taimi\c{n}a.
\newblock \emph{Crocheting adventures with hyperbolic planes}.
\newblock CRC Press, Boca Raton, FL, New York, second edition, 2018.
\newblock ISBN 978-1-138-30115-3.

\bibitem[Thurston(2022)]{Thurston2022}
William~P Thurston.
\newblock \emph{The Geometry and Topology of Three-Manifolds: With a Preface by
  Steven P. Kerckhoff}, volume~27.
\newblock American Mathematical Society, 2022.

\bibitem[Turner and van~de Griend(1996)]{Turner1996}
J.~C. Turner and P.~van~de Griend, editors.
\newblock \emph{History and science of knots}, volume~11 of \emph{Series on
  Knots and Everything}.
\newblock World Scientific Publishing Co., Inc., River Edge, NJ, 1996.
\newblock ISBN 981-02-2469-9.
\newblock \doi{10.1142/9789812796134}.
\newblock URL \url{https://doi.org/10.1142/9789812796134}.

\bibitem[van Wijk and Cohen(2016)]{Wijk2016}
Jarke~J. van Wijk and Arjeh~M. Cohen.
\newblock Visualization of {S}eifert surfaces.
\newblock In \emph{Six papers on signatures, braids and {S}eifert surfaces},
  volume~30 of \emph{Ensaios Mat.}, pages 217--245. Soc. Brasil. Mat., Rio de
  Janeiro, 2016.
\newblock Reprinted from IEEE Transactions on Visualization and Computer
  Graphics {{\bf{1}}2} (2006), no. 4, 485--496.

\end{thebibliography}

\end{document}